\documentclass[a4paper]{scrartcl}

\usepackage[utf8]{inputenc}
\usepackage[T1]{fontenc}
\usepackage{lmodern}
\usepackage{amsmath, amssymb}
\usepackage{amsthm}
\usepackage{graphicx}
\usepackage{float}
\usepackage{caption}
\usepackage{subcaption}
\usepackage[backend=bibtex,style=authoryear,natbib]{biblatex}
\usepackage[dvipsnames,svgnames]{xcolor}
\usepackage{tikz,pgf}
\usepackage{dirtree}
\usetikzlibrary{calc}
\usetikzlibrary{mindmap}
\usetikzlibrary{shadows}
\usetikzlibrary{fadings}
\usetikzlibrary{backgrounds}
\usetikzlibrary{arrows.meta}
\usetikzlibrary{svg.path}
\usepackage{xurl}

\usepackage{booktabs}
\usepackage{listings}
\usepackage[colorlinks=true,linkcolor=pyrigiblue,citecolor=PineGreen,urlcolor=pyrigiblue]{hyperref}
\usepackage[nameinlink]{cleveref}

\colorlet{colbg}{white}
\colorlet{colfg}{black}
\colorlet{colgraphv}{colfg!75!colbg}
\colorlet{colgraphe}{colfg!55!colbg}
\colorlet{colG}{DarkSeaGreen}
\definecolor{colR}{HTML}{CC6677}
\definecolor{colO}{HTML}{DDCC77}
\definecolor{colB}{HTML}{6699CC}
\colorlet{colY}{Gold!70!black}
\colorlet{colGray}{white!60!black}
\definecolor{pyrigiblue}{HTML}{4584b6}
\definecolor{pyrigiyellow}{HTML}{ffde57}
\definecolor{pyrigigray}{HTML}{646464}

\colorlet{darkGreen}{colG!70!black}
\colorlet{darkBlue}{colB!70!black}
\colorlet{darkRed}{colR!70!black}
\colorlet{darkGray}{black!70!white}

\newcommand\pythonstyle{\lstset{
		language=Python,
		basicstyle=\ttfamily,
		morekeywords={self},
		keywordstyle=\bfseries\ttfamily\color{darkBlue},
		emph={},
		emphstyle=\ttfamily\bfseries\color{darkRed},
		stringstyle=\ttfamily\color{darkGreen},
		commentstyle=\ttfamily\color{darkGray},
		frame=tb,
		showstringspaces=false,
		columns=flexible
}}
\lstnewenvironment{python}[1][frame=none,xleftmargin=.25in]
{
	\pythonstyle
	\lstset{#1}
}
{}

\lstnewenvironment{pythonfig}[1][]
{
	\pythonstyle
	\lstset{#1}
}
{}

\newcommand\pythoninline[1]{{\pythonstyle\lstinline!#1!}}

\theoremstyle{definition}

\tikzstyle{vertex}=[fill=colgraphv,circle,inner sep=0pt, minimum size=4pt]
\tikzstyle{edge}=[line width=1.5pt,colgraphe]
\tikzstyle{labelsty}=[font=\scriptsize]
\tikzstyle{indicatededge}=[pin={[pin distance=6pt,pin edge={thin,path fading=east,fading transform={rotate=#1}}]20:},pin={[pin distance=6pt,pin edge={thin,path fading=east,fading transform={rotate=#1}}]-20:},pin={[pin distance=10pt,pin edge={thin,path fading=east,fading transform={rotate=#1}}]2:}]
\tikzstyle{construle}=[ultra thick,-{Classical TikZ Rightarrow[]}]
\tikzstyle{edgeq}=[edge,dashed]

\tikzstyle{genericgraph}=[dashed,black!70!white]

\tikzstyle{elem}=[draw=black!80!white,rounded corners,inner sep=0.3cm,line width=1pt]
\tikzstyle{head}=[anchor=west,fill=white,inner sep=2pt]
\tikzstyle{ref}=[line width=1.5pt, dashed,-{LaTeX[round]},colB,rounded corners]
\tikzstyle{sref}=[ref,line width=0.5pt,colB!50!white]
\tikzstyle{refrev}=[line width=1.5pt, dashed,{LaTeX[round]}-,colB,rounded corners]
\tikzstyle{srefrev}=[refrev,line width=0.5pt,colB!50!white]
\tikzstyle{con}=[line width=1.5pt,-{LaTeX[round]},colR,rounded corners]
\tikzstyle{conb}=[line width=2.5pt,white,rounded corners]

\tikzstyle{frame}=[line width=1.5pt,black!50!white,rounded corners]
\tikzstyle{frametext}=[black!50!white,fill=white,anchor=west,font=\small]

\newcommand{\pyrigilogo}[1]{%
	\pgfmathparse{#1*2.09414}\let\lw=\pgfmathresult
	\pgfmathparse{#1*2.09414/2}\let\sl=\pgfmathresult
	\begin{tikzpicture}[scale=#1,yscale=-1,tpedge/.style={pyrigigray,line width=\lw pt,shorten <= -\sl pt,shorten >= -\sl pt,{Round Cap[]}-{Round Cap[]}}]
		\draw[tpedge] svg {M 98.79062,120.93986 84.022775,134.60024};
		\draw[tpedge] svg {m 98.794581,120.93194 18.841409,20.41147};
		\draw[tpedge] svg {m 117.63897,141.37349 -14.79746,13.74874};
		\draw[tpedge] svg {m 84.089545,134.70916 18.673085,20.39718};

		\fill[pyrigiblue] svg {m 92.802948,93.272202 c 0.322655,1.068806 0.609308,1.860968 0.776946,2.917297 0.226705,1.428535 0.279535,1.969444 0.262276,3.605163 -0.02922,2.769098 -0.585757,4.768378 -1.670333,6.556428 -1.612422,2.65825 -3.556782,5.00722 -5.565609,7.38016 -2.183482,2.57926 -4.848813,4.92153 -6.810334,7.67336 -2.77676,3.89553 -3.911539,6.52585 -4.78382,9.92378 -0.572888,2.23166 -0.687265,4.29225 -0.625978,6.68309 0.05551,2.16531 0.169428,3.66093 1.047397,6.37299 0.650901,2.01065 1.629267,3.88444 2.737763,5.68379 1.330881,2.16033 3.07493,3.89816 5.000805,5.55021 2.020002,1.73281 4.362691,3.2641 6.825757,4.27221 1.465424,0.59978 3.223195,1.29141 5.007728,1.60867 2.207325,0.39242 3.613246,0.48177 5.187974,0.46636 0.0798,-0.003 0.035,-0.004 0.076,5.8e-4 0.17671,-0.10751 0.58809,-0.44747 0.0989,-0.52004 -0.489201,-0.0726 -1.923573,-0.12571 -2.257509,-0.17436 -1.638796,-0.23889 -3.380034,-0.67884 -4.907031,-1.31994 -1.398642,-0.58722 -2.739783,-1.31516 -3.978306,-2.191 -1.393559,-0.98546 -2.653746,-2.15609 -3.760433,-3.45547 -1.044645,-1.22653 -1.953616,-2.56856 -2.705022,-3.99371 -0.688645,-1.2816 -1.188154,-2.66467 -1.477324,-4.09054 -0.09974,-0.49183 -0.174544,-0.98895 -0.224161,-1.48833 -0.136574,-0.92791 -0.186512,-1.86847 -0.149022,-2.80564 0.02602,-0.65034 0.09405,-1.29831 0.177003,-1.94387 0.0984,-0.76584 0.217806,-1.52959 0.386386,-2.28311 0.277256,-1.23927 0.686431,-2.44697 1.17991,-3.61707 0.952804,-2.25922 2.220526,-4.37942 3.69001,-6.34221 0.729222,-0.97402 1.507768,-1.9097 2.286881,-2.8443 1.778585,-2.13353 3.561891,-4.26417 5.266056,-6.4576 1.184911,-1.5251 2.332912,-3.08243 3.344176,-4.72781 0.678756,-1.10437 1.286198,-2.39822 1.685167,-3.62203 l 0.03869,-0.0144 c -0.184518,-0.14757 0.547356,0.61044 0.694006,0.79522 0.460063,0.57972 0.892313,1.19506 1.305093,1.80934 0.5899,0.87786 0.97639,1.74017 1.4168,2.70175 0.3171,0.69235 0.60141,1.47194 0.81882,2.20176 0.19644,0.65941 0.31107,1.15122 0.39283,1.8344 0.0491,0.41071 0.11794,1.08128 0.12646,1.7912 0.0103,0.85748 -0.0238,1.2982 -0.0137,2.16264 0.0103,0.8801 0.10528,1.84523 0.23451,2.71585 0.19451,1.31068 0.53136,2.83296 1.20965,3.97123 1.01378,1.70128 2.67955,3.04446 4.60258,3.51782 0.51069,0.12572 1.27131,0.15459 1.79721,0.16045 0.78728,-0.074 1.40617,-0.27575 2.13601,-0.50034 0.66883,-0.18908 1.30125,-0.55915 1.90028,-0.91162 0.48456,-0.28511 0.93954,-0.62038 1.35503,-0.99913 0.25794,-0.22844 0.49526,-0.47993 0.70863,-0.75046 0.10547,-0.13371 0.20507,-0.2722 0.29837,-0.41467 0.47765,1.09656 0.8538,2.23727 1.12174,3.40296 0.14172,0.61654 0.25334,1.24007 0.3342,1.86748 0.35208,1.83208 0.46126,3.7108 0.32382,5.57133 -0.0855,1.15758 -0.26615,2.30723 -0.51906,3.44007 -0.28642,1.28296 -0.66618,2.54603 -1.15318,3.76703 -0.78979,1.98017 -1.90259,3.87996 -3.21516,5.55984 -1.05657,1.2173 -1.76708,2.14327 -2.99024,3.14938 -0.92726,0.76272 -2.00309,1.50241 -3.04861,2.0957 -1.37024,0.77756 -2.20918,1.00089 -3.74647,1.50039 -1.17613,0.38215 -2.74538,0.72901 -3.89776,0.75937 -1.0727,0.0283 -0.96689,0.59814 -0.88851,0.69218 1.8059,-0.009 3.6267,-0.24607 6.36219,-0.90181 2.26746,-0.54354 4.60731,-1.61674 6.62143,-2.79154 3.25881,-1.90081 5.67136,-4.50685 7.71126,-7.47567 1.51157,-2.19989 2.92888,-6.07426 3.39541,-8.55582 0.97695,-5.19662 0.83641,-9.96011 0.023,-14.00091 -0.52595,-2.6128 -1.29952,-5.02033 -2.42886,-7.37314 -1.03519,-2.15669 -2.57189,-4.56723 -4.75628,-6.45939 0.2256,0.7398 0.34591,1.56712 0.39845,2.05139 0.0509,0.46865 0.12101,1.17121 0.12731,1.74782 0.009,0.82284 -0.0129,1.49391 -0.0763,2.34172 -0.064,0.85483 -0.23161,1.62023 -0.40474,2.42906 -0.13256,0.61923 -0.35616,1.3659 -0.64456,1.92969 -0.3326,0.6502 -0.79456,1.22327 -1.36266,1.73015 -0.50044,0.4465 -1.26685,0.86494 -2.08033,1.02887 -0.64345,0.12966 -1.35809,0.16373 -2.00388,0.0463 -0.72165,-0.21094 -1.3147,-0.67892 -1.81557,-1.22985 -0.56639,-0.623 -0.79135,-1.23361 -1.00627,-2.04768 -0.28268,-1.07071 -0.32711,-2.16221 -0.33558,-3.26958 0.009,-0.99498 0.0122,-1.99042 0.0187,-3.02523 0.0378,-0.58876 -0.005,-1.03601 0.006,-1.62589 0.0147,-0.78212 -0.0708,-1.75506 -0.13832,-2.53439 -0.16937,-1.95185 -0.67252,-3.86682 -1.39338,-5.68859 -0.52652,-1.33057 -1.16883,-2.61379 -1.88352,-3.85349 -0.81373,-1.41148 -1.64534,-2.70369 -2.70453,-3.94165 -1.03368,-1.208146 -2.48576,-2.734448 -4.183216,-4.086091 -1.660572,-1.322273 -3.656156,-2.325598 -4.874781,-3.116943 z};

		\fill[pyrigiyellow] svg {m 101.49864,128.05639 c -0.47969,0.36484 -1.30384,1.05485 -1.695524,1.48793 -0.592857,0.65553 -1.202494,1.58174 -1.64251,2.34828 -0.561111,0.97751 -0.859378,2.10851 -1.141933,3.19963 -0.26941,1.04036 -0.425891,2.11231 -0.426944,3.18699 -4.73e-4,0.48484 0.03072,0.96949 0.08552,1.45122 0.05838,0.51319 0.143433,1.02325 0.199006,1.53675 0.07907,0.73064 0.09682,1.47899 -0.0852,2.19099 -0.09742,0.38105 -0.254623,0.7526 -0.505228,1.05572 -0.406069,0.49119 -1.035445,0.76639 -1.670044,0.82498 -0.340116,0.0314 -0.6869,0.005 -1.013912,-0.0936 -0.327013,-0.0986 -0.633782,-0.27061 -0.877464,-0.50994 -0.37995,-0.37318 -0.779521,-1.15181 -0.860156,-1.67822 -0.09428,-0.61545 -0.08796,-1.50689 0.03847,-2.4208 0.126423,-0.9139 0.300842,-1.91802 0.42301,-2.54474 -0.473175,0.28446 -1.624315,1.79419 -1.919507,2.35221 -0.309813,0.58566 -0.829613,1.42922 -1.072381,2.0457 -0.331547,0.84192 -0.641539,1.71291 -0.829009,2.61941 -0.366374,1.77157 -0.366743,4.2968 0.0055,6.40219 0.283801,1.60527 0.8473,3.16795 1.717519,4.54641 0.824855,1.30662 1.919607,2.43748 3.17692,3.33569 1.226771,0.87639 2.616814,1.53615 4.089291,1.85995 0.932339,0.20501 2.118796,0.29161 3.072486,0.24971 1.09212,-0.048 1.99859,-0.19998 3.04047,-0.53093 1.07838,-0.34256 2.0628,-0.8906 3.01098,-1.50799 0.83897,-0.59808 1.60161,-1.3031 2.26364,-2.09258 1.20244,-1.434 2.07296,-3.16481 2.38433,-5.01014 0.14349,-0.85036 0.16868,-1.71773 0.12419,-2.57895 -0.0774,-1.49924 -0.27312,-2.96296 -0.80156,-4.36811 -0.51476,-1.3688 -1.22788,-2.72194 -2.11059,-3.88789 -0.50843,-0.6716 -1.0803,-1.32312 -1.62609,-1.96475 -0.78805,-0.92641 -1.51718,-1.5832 -2.25469,-2.55034 -0.68444,-0.89755 -1.31851,-2.00897 -1.7137,-3.06629 -0.41402,-1.10767 -0.67137,-2.32576 -0.75692,-3.42963 -0.0633,-0.8162 -0.0468,-1.85222 0.14827,-2.45152 0.26662,-0.81906 -0.10117,-0.48792 -0.49804,-0.20538 z m -0.48219,16.72507 c 0.78277,0.737 1.51102,1.53171 2.17684,2.37585 0.87958,1.11512 1.65167,2.31908 2.2476,3.60826 0.55711,1.20521 0.87687,2.49427 0.96938,3.81879 0.0316,0.45232 0.0668,0.8902 0.0277,1.34195 -0.0563,0.65083 -0.19874,1.31262 -0.51188,1.88596 -0.32806,0.60065 -0.59587,0.91856 -1.0441,1.4053 -0.31935,0.3468 -0.7145,0.63444 -1.10802,0.89406 -0.35444,0.23384 -0.92527,0.4603 -1.45738,0.58879 -0.84116,0.20311 -1.34921,0.2644 -2.04709,0.25013 -0.838021,-0.0171 -1.276703,-0.0118 -2.232264,-0.26993 -0.572697,-0.15474 -0.931132,-0.2813 -1.43395,-0.56118 -0.441566,-0.24578 -0.685146,-0.56866 -0.940082,-0.84755 -0.576777,-0.63097 -0.891199,-1.41063 -1.05754,-1.98112 -0.130711,-0.44829 -0.156606,-0.92142 -0.160968,-1.38836 -0.0056,-0.57459 0.0496,-1.15512 0.231395,-1.70022 0.177067,-0.53094 0.470945,-1.01725 0.8149,-1.45878 0.506,-0.64953 1.120074,-1.2065 1.750528,-1.73606 0.425822,-0.35768 0.860758,-0.70485 1.27046,-1.08088 0.251558,-0.23089 0.49329,-0.4724 0.728136,-0.72027 0.266194,-0.28096 0.524387,-0.57094 0.744784,-0.8891 0.248821,-0.35918 0.447731,-0.75208 0.604891,-1.15979 0.0869,-0.22545 0.16115,-0.45584 0.22039,-0.69009 0.0332,-0.13135 0.0617,-0.26384 0.0862,-0.39707 0.0451,-0.24507 0.0769,-0.49251 0.0952,-0.74103 z};

		\draw[tpedge] svg {m 84.357142,134.89207 12.93376,-0.7066 1.733779,-12.75578};
		\draw[tpedge] svg {m 102.7593,154.94978 1.38156,-13.37205 13.0373,-0.41673};
		\draw[tpedge] svg {m 97.510769,134.36601 6.609361,7.17695};
	\end{tikzpicture}
}

\newcommand{\pyrigi}{\textsc{PyRigi}}

\newcommand{\RR}{\mathbb{R}}

\title{PyRigi --- a general-purpose Python package for the rigidity and flexibility of bar-and-joint frameworks}
\author{%
	Matteo Gallet%
	\thanks{University of Trieste, Department of Mathematics, Informatics and Geosciences}
	\and
	Georg Grasegger%
	\thanks{Johann Radon Institute for Computational and Applied Mathematics (RICAM), Austrian Academy of Sciences}
	\and
	Matthias Himmelmann%
	\thanks{Institute for Computational and Experimental Research in Mathematics (ICERM), Brown University in Providence, RI, USA}
	\and
	Jan Legerský%
	\thanks{Department of Applied Mathematics, Faculty of Information Technology, Czech Technical University in Prague}
	}
\date{}

\begin{document}
\maketitle

\begin{abstract}
  We present \pyrigi{}, a novel Python package designed to study the rigidity properties of graphs and frameworks.
  Among many other capabilities, \pyrigi{} can determine whether a graph admits only finitely many ways, up to isometries, of being drawn in the plane once the edge lengths are fixed,
  whether it has a unique embedding,
  or whether it satisfies such properties even after the removal of any of its edges.
  By implementing algorithms from the scientific literature, \pyrigi{} enables the exploration of rigidity properties of structures that would be out of reach for computations by hand.
  With reliable and robust algorithms, as well as clear, well-documented methods that are closely connected to the underlying mathematical definitions and results, \pyrigi{} aims to be a practical and powerful general-purpose tool for the working mathematician interested in rigidity theory.
  \pyrigi{} is open source and easy to use, and awaits researchers to benefit from its computational potential.
\end{abstract}

\begin{center}
	\pyrigilogo{1} \tikz{\node[align=center] {{\color{pyrigigray}\huge\pyrigi} \\[2pt]%
			\href{https://pyrigi.github.io}{\hspace{0.5pt}\foreach \t in {p,y,r,i,g,i,.,g,i,t,h,u,b,.,i,o} {\t\kern0.175pt}}}}
\end{center}

\section{Introduction}

This paper serves as an introduction to \pyrigi{},
a general-purpose Python package. Its aim is to help researchers in rigidity and flexibility theory in their everyday work
by implementing various algorithmic and combinatorial aspects of the mathematical theory, allowing users to perform computations that would be infeasible to carry out by hand.
This document presents an overview of \pyrigi{}’s functionality, design philosophy, and implementation details.

Rigidity theory is the study of graph-theoretic, matroid-theoretic, geometric, and metric properties related to the various forms of structural rigidity.
Concretely, one starts from a simple, loopless graph and assigns to each of its vertices a point in $\mathbb{R}^d$, namely a \emph{realization};
in this way, a \emph{framework} is obtained.
Then, one can ask whether it is possible to assign other realizations to the same graph
so that the Euclidean distance between points corresponding to an edge of the graph is the same as in the original realization.
Clearly, this is always possible by applying an isometry, i.e., a rototranslation, to the original realization.
If this is the only possible way to obtain such ``compatible'' realizations in a sufficiently small neighborhood,
then the framework is called \emph{rigid}; it is called \emph{flexible} otherwise.
In addition to this first notion of rigidity, many other can be introduced:
for example, if a framework is rigid, does there exist a unique realization (up to isometries) with given distances between points connected by an edge, or are there finitely many of them?
In the first case, the framework is called \emph{globally rigid}.
If the framework is rigid and removing any edge preserves this property, then it is called \emph{redundantly rigid}.
All these notions can be transferred from frameworks to abstract graphs by considering what the behavior of a ``generic'' framework is.

Rigidity theory sits at the interface between graph theory, matroid theory, and metric geometry,
and benefits from contributions from algebraic and tropical geometry, functional analysis,
semidefinite programming, and symbolic computation.
In addition to its theoretical interest, rigidity theory plays a role in applications,
such as protein and molecule analysis \parencite{ThorpeDuxbury2002,HermansPflegerEtAl2017}, formation control \parencite{ZelazoZhao2019},
matrix completion \parencite{SingerCucuringu2010}, foldable structures \parencite{Tachi2016}, chemistry and modern art \parencite{Graver2001}.

The area of rigidity theory is currently experiencing an exciting period of growth.
For example, the average number of yearly entries tagged with MSC code 52C25
``Rigidity and flexibility of structures (aspects of discrete geometry)'' on
the \emph{zbMATH Open} abstracting and reviewing service\footnote{\url{https://zbmath.org}}
increased from 28.4 in the five-year period 2015--2019 to 50.6 in 2019--2024.
Moreover, several major international scientific gatherings focusing on rigidity
theory and related topics have been organized in recent years.
These include the ``Thematic Program on Geometric Constraint Systems,
Framework Rigidity, and Distance Geometry'' (2021)
and the subsequent ``Focus Program on Geometric Constraint Systems'' (2023),
both hosted by the Fields Institute for Research in Mathematical Sciences
(Toronto, Canada), the ``Special Semester on Rigidity and Flexibility'' (2024)
at RICAM (Linz, Austria) and the program ``Geometry of Materials, Packings,
and Rigid Frameworks'' (2025) at ICERM (Providence, USA).

Parallel to the scientific activity, the last 15 years have seen the emergence of a wide range of software tools devoted to rigidity theory. Their scope ranges from projects taylored to applications
(as for example \emph{FIRST/MSU ProFlex} \parencite{ProFlexPaper}, \emph{KINARI} \parencite{KinariPaper, KinariWebPaper}, or \textcite{Algorithm951}),
to code developed as supporting material for research papers or scientific work
\parencite{Anand2024,Bartzos2020,RealizationCountSoftware,Dance2024,DewarSoftware,SphereCountSoftware,CalligraphSoftware,FlexrilogPaper,HaghSadjadi2022,Islam2024,McGlue2020,McKenzie2020,Algorithm990,Poetzsch2024,GPART}.

Despite the growing interest in developing software tools for rigidity theory,
it seemed to us that there is a lack of general-purpose software to support researchers in their everyday work. Such a software package should contain features that allow testing for various rigidity properties,
be able to manipulate graphs with rigidity-oriented operations such as $k$-extensions or coning,
and offer databases with graphs that have interesting rigidity properties available.

This motivation led to the workshop ``Code of Rigidity'' during the 2024 Special
Semester on Rigidity and Flexibility in Linz which was dedicated to presenting and
spreading knowledge about software in rigidity theory.
A survey was prepared and circulated before the workshop among possibly interested researchers
about their needs and wishes regarding software for rigidity theory and the answers
were the starting point for the workshop discussions.
The (approx.\ 25) workshop participants agreed on the importance and usefulness of starting to develop a general-purpose software platform for rigidity theory and set associated goals.
Via rounds of consultations, it was decided that the software should be standalone (and not embedded as a package of a wider already existing project) and developed in Python (which accumulated the most potential users and contributors).
The name \pyrigi{} was selected from a shortlist of proposals through a vote by the workshop participants as a portmanteau combining ``Python'' and ``Rigidity''.
Finally, the four authors of this paper volunteered to act as maintainers of the project.

After a year, \pyrigi{} has reached its first stable release, and this article is meant to introduce the readers to its features and (current) limitations,
and to hopefully convince researchers in rigidity theory or nearby fields to consider adopting it in their work.
The main strengths of \pyrigi{} can be summerized as follows:
\begin{itemize}
	\item it is \emph{general-purpose}: a user finds a wealth of commands, all in the same software;
	\item it is \emph{easy to use}: as a Python package, it is easy to install on many platforms (simply by running \texttt{pip install pyrigi}) and it provides an intuitive interface, with method names derived from the mathematical concepts they implement;
	\item it is the outcome of \emph{community involvement}: it is developed as a \emph{GitHub} repository\footnote{\url{https://github.com/PyRigi/PyRigi}}, so that anyone can contribute, and it has a Zulip channel\footnote{\url{https://pyrigi.zulipchat.com}} to discuss issues, feature requests, and ongoing development; it originated from a workshop, and several events have already taken place to enhance it;
	\item it is \emph{well-documented}, including \emph{mathematical background}: we put an emphasis on providing detailed interface instructions and to link them to the related mathematical concepts;
	\item it \emph{implements scientific results} from the recent literature.
\end{itemize}

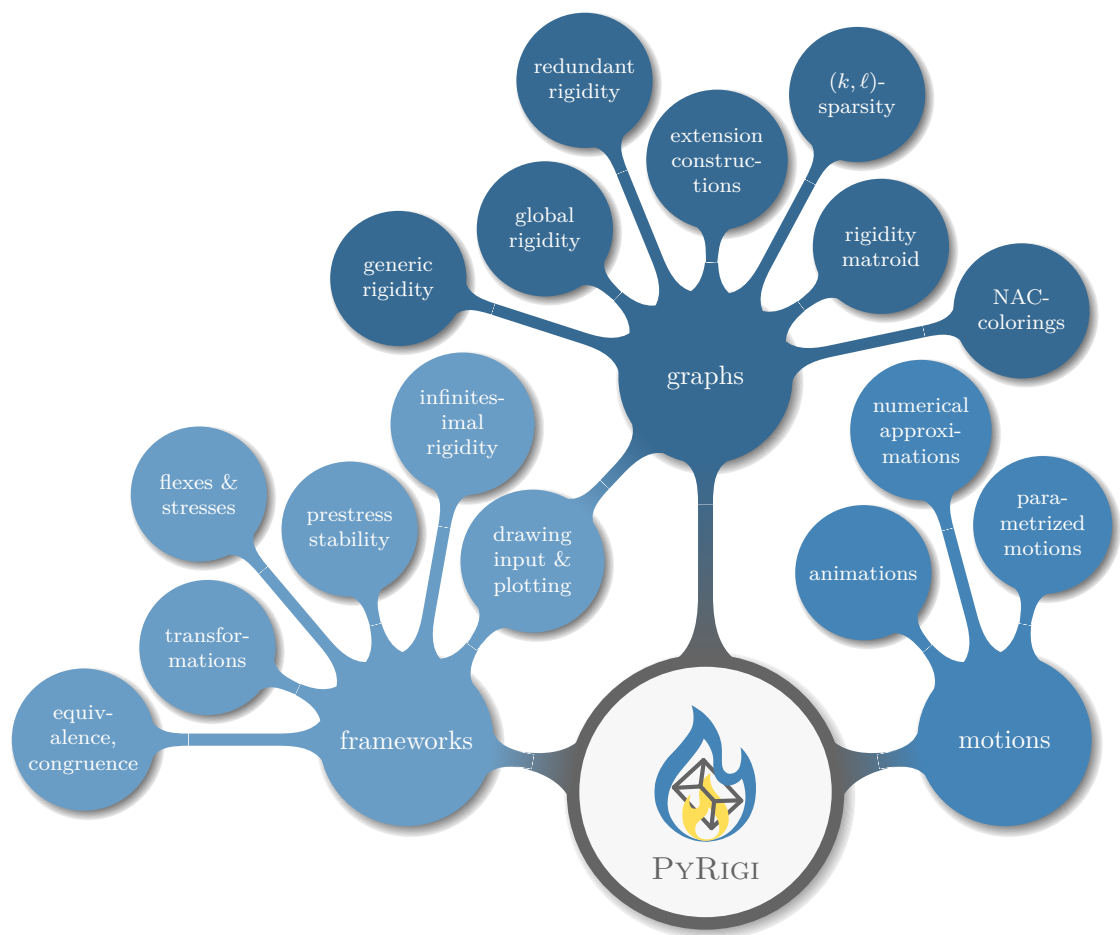
\begin{figure}[ht]
	\centering
	\begin{tikzpicture}[mindmap,scale=0.85]
		\begin{scope}[
			every node/.style={concept, circular drop shadow, execute at begin node=\hskip0pt},
			root concept/.append style={concept color=pyrigigray, fill=pyrigigray!5!white, line width=1ex, text=pyrigigray,minimum size=2cm,align=center,inner sep=0pt},
			text=white,
			concept color=pyrigigray,
			graph/.style={concept color=pyrigiblue!80!black},
			framework/.style={concept color=pyrigiblue!80!white},
			motion/.style={concept color=pyrigiblue},
			grow cyclic,
			level 1/.append style={level distance=4.7cm,sibling angle=80},
			level 2/.append style={level distance=3.4cm,sibling angle=25, font=\scriptsize}]
			]
			\node[root concept] (rt) {}
			[clockwise from=170]
			child[framework]{ node {frameworks}
				[clockwise from=180]
				child[level distance=5cm]{ node {equivalence, congruence}}
				child{ node {transformations}}
				child[level distance=5cm]{ node {flexes \& stresses}}
				child{ node {prestress stability}}
				child[level distance=5cm]{ node (ir) {infinitesimal rigidity}}
				child{ node (gd) {drawing input \& plotting}}
			}
			child[graph]{ node[yshift=1.5cm] (g) {graphs}
				[clockwise from=162]
				child[level distance=5cm]{ node (gr) {generic rigidity}}
				child{ node {global rigidity}}
				child[level distance=5cm]{ node {redundant rigidity}}
				child{ node {extension constructions}}
				child[level distance=5cm]{ node {$(k,\ell)$-sparsity}}
				child{ node {rigidity matroid}}
				child[level distance=5cm]{ node {NAC-colorings}}
			}
			child[motion]{ node {motions}
				[counterclockwise from=80]
				child{ node {para\-metrized motions}}
				child[level distance=5cm]{ node {numerical approximations}}
				child{ node {animations}}
			};
			\path (g) to[circle connection bar switch color=from (pyrigiblue!80!black) to (pyrigiblue!80!white)] (gd);
		\end{scope}
		\node[pyrigigray] at (rt) {\pyrigilogo{0.75}\\\pyrigi};
	\end{tikzpicture}
	\caption{Overview of \pyrigi's functionality}
	\label{fig:pyrigi-overview}
\end{figure}

\Cref{fig:pyrigi-overview} provides a glimpse of the current capabilities of \pyrigi{}.
As one may notice, \pyrigi{} rests on the three pillars of graphs, frameworks, and motions.
Each of these pillars corresponds to a Python class, that exposes methods
which can be used to create and modify graphs and frameworks,
both by explicitly entering their vertices, edges, and realizations,
and by fetching them from databases of known examples;
users are able to call methods to check for infinitesimal and generic rigidity, in addition to their global and redundant counterparts.
A non-exhaustive list of possible constructions on graphs includes $k$-extensions, intersections, sums, and coning.
Extensive plotting features are provided, and one can numerically generate motions of flexible frameworks.

\pyrigi{} aims at being easy to use, especially for people with little to no experience with symbolic computations and programming:
our chosen programming language, Python, is known to be easy to learn also for beginners,
and we strive to provide clear and detailed documentation for the introduced objects.
In addition, the package documentation is enhanced by interlacing it with mathematical definitions,
that provide precise mathematical meaning to the inputs and outputs of the Python computations.

\pyrigi{} is an open-source software package.
We believe that the outcomes of research should be publicly and freely accessible,
in accordance to the FAIR Principles to Research Software \parencite{Barker2022},
that is why all \pyrigi{}'s code is accessible via its GitHub page\footnote{\url{https://github.com/PyRigi/PyRigi}}
and is stored in Zenodo \parencite{zenodo}.
Everyone is welcome to join this project.
\pyrigi{} is released under an MIT license, which grants several rights to use, modify and copy the software.
We kindly ask users and developers to carefully read the license and be aware of its consequences.

This article serves as an introduction to \pyrigi{} by describing its core functionalities, guiding philosophy, and internal architecture.
\Cref{code_math} highlights \pyrigi{}'s unique feature of integrating formal mathematical statements with its interface documentation to give a precise meaning to the quantities used in its methods.
\Cref{functionality} provides a detailed overview of the functionalities implemented in \pyrigi{}.
\Cref{implementation} sheds light on the implementation details of the main classes and methods in \pyrigi{} and their internal structure.
Finally, \Cref{conclusions} outlines potential directions for future development and extensions of \pyrigi{}.

\section{Code meets Math}
\label{code_math}

\begin{figure}[ht]
	\centering
	\begin{tikzpicture}[]
		\node[elem,anchor=north west] (doc) at (0,0) {\includegraphics[width=5cm]{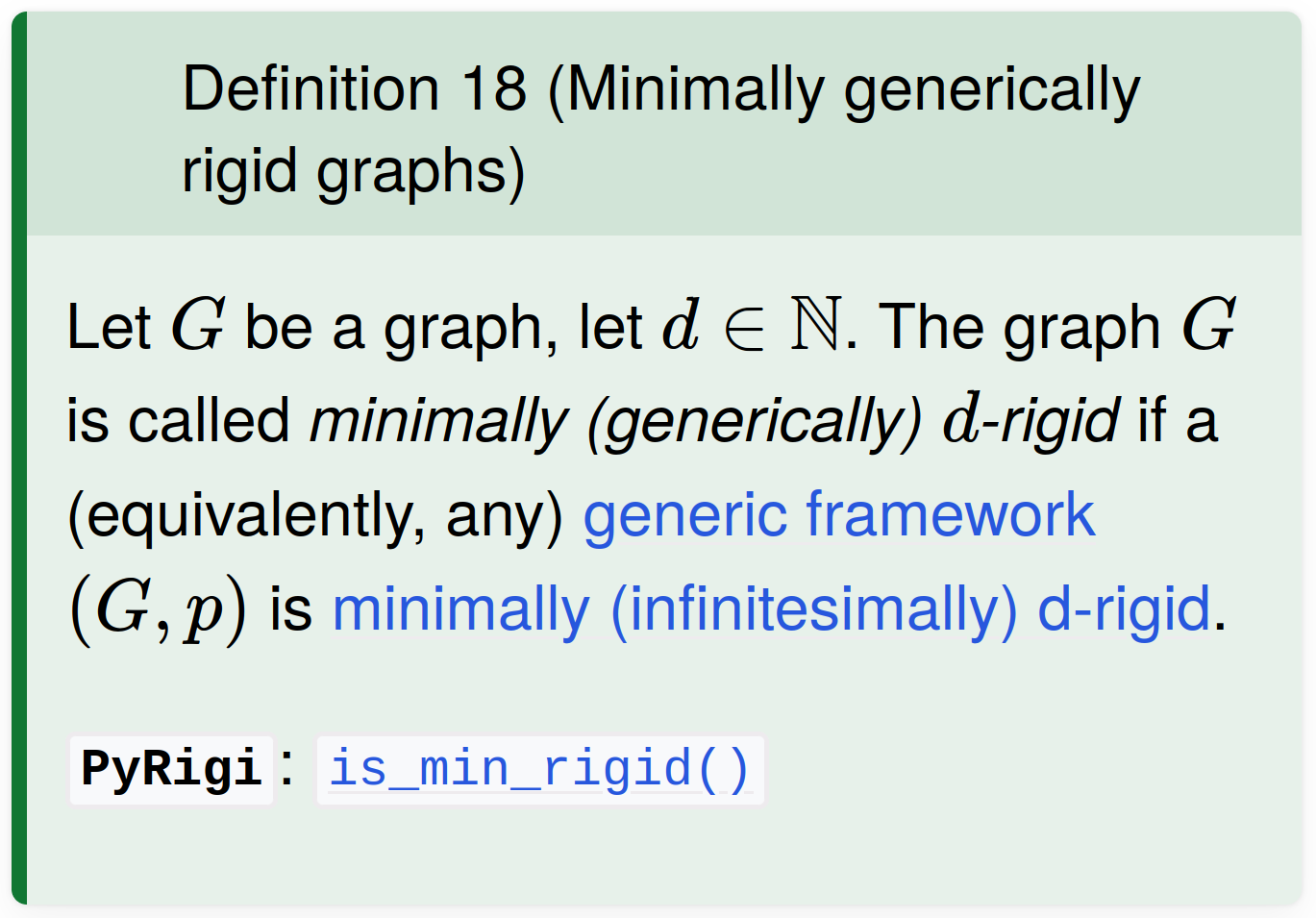}};
		\node[elem,anchor=north east] (thm) at (13.5,0) {\includegraphics[width=5cm]{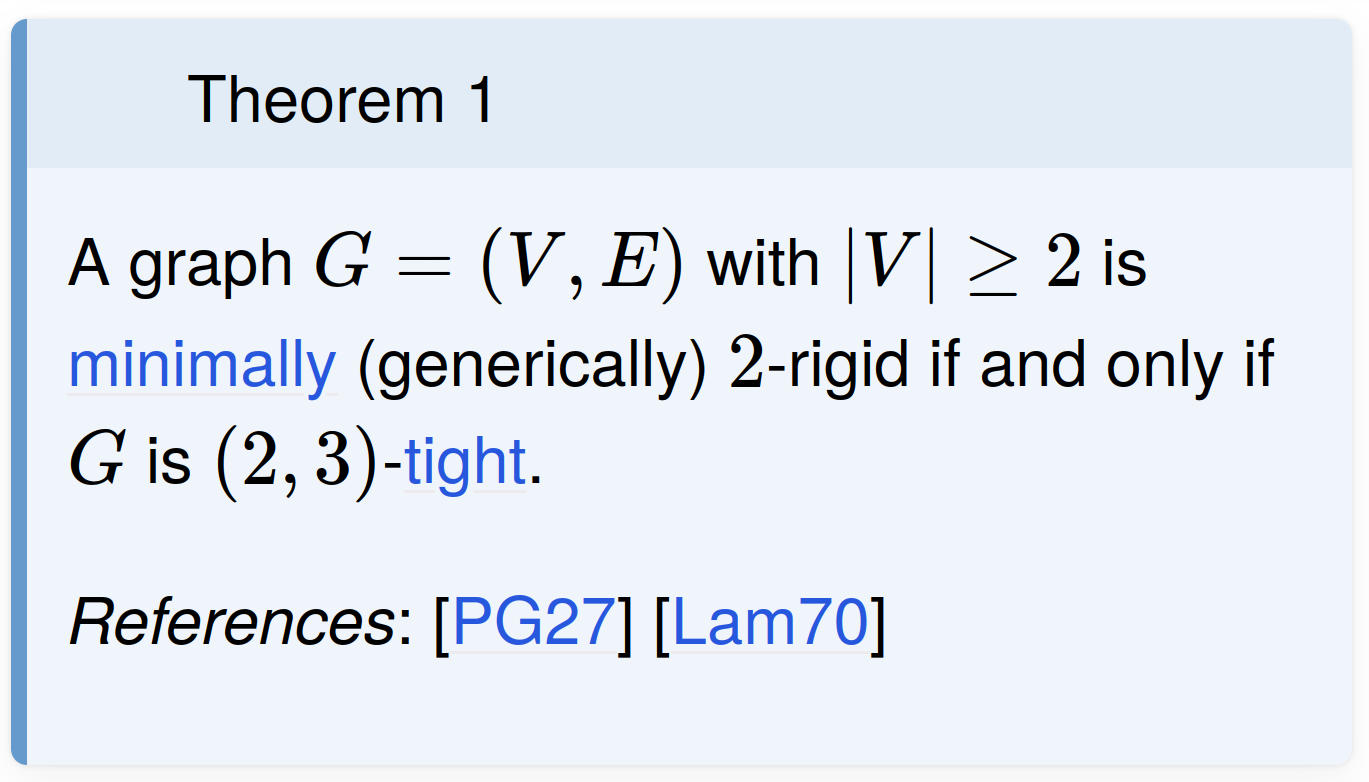}};
		\node[elem,anchor=north east,minimum width=5cm,minimum height=1.5cm] (cdoc) at (13.5,-4) {\includegraphics[width=4cm,clip=true,trim=0cm 25cm 0cm 0cm]{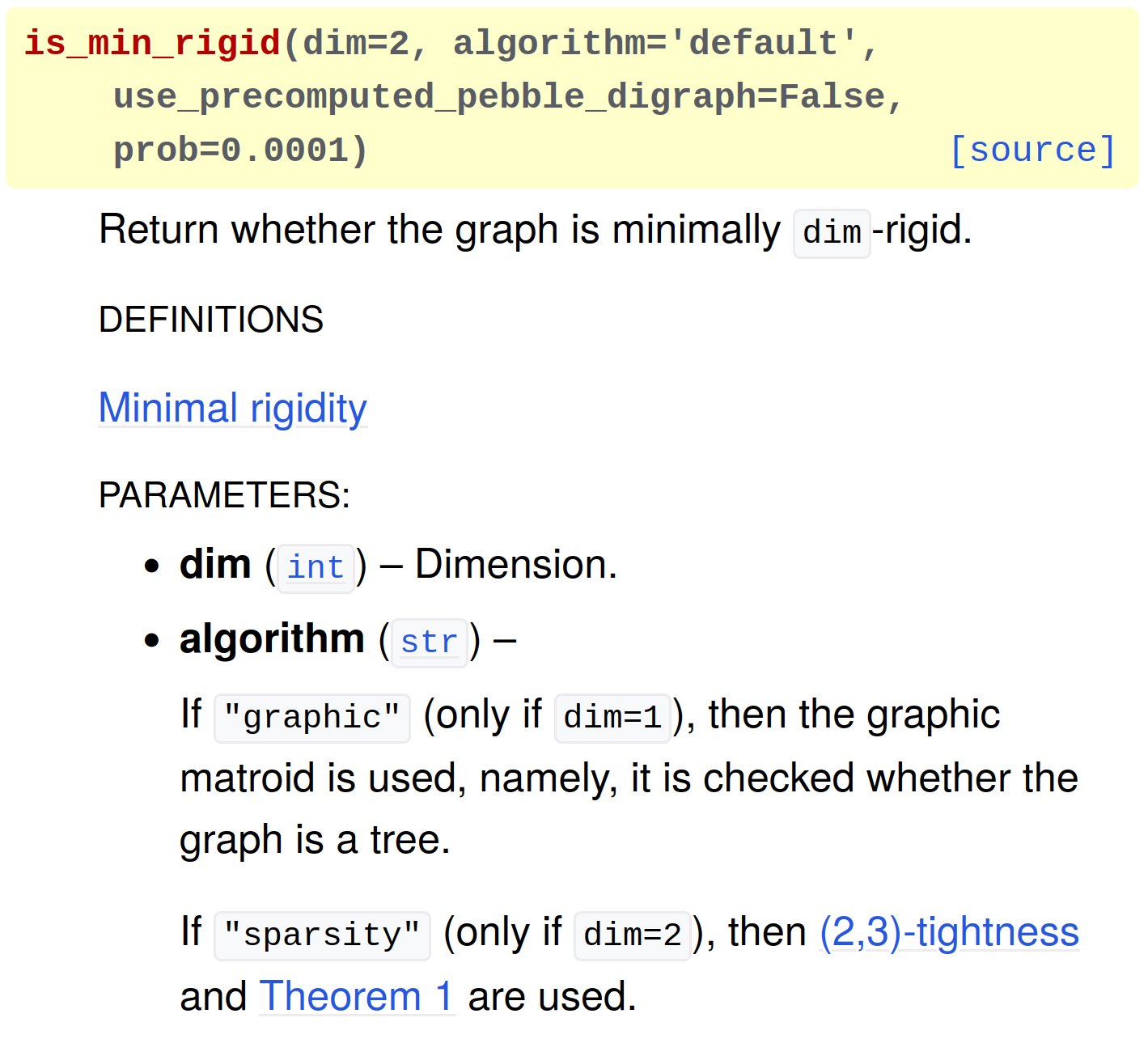}};
		\node[elem,anchor=north east] (ex) at (13.5,-6.25) {%
			\begin{tikzpicture}[%
				scale=1.2,
				gvertex/.style={fill=black,draw=white,circle,inner sep=0pt,minimum size=4pt},
				edge/.style={line width=1.5pt,black!60!white}
			]
				\node[gvertex] (0) at (0.12262, -0.16979) {};
				\node[gvertex] (1) at (0.01795, -1.0) {};
				\node[gvertex] (2) at (0.94358, -0.33192) {};
				\node[gvertex] (3) at (-0.60737, 0.84134) {};
				\node[gvertex] (4) at (-0.75238, -0.04083) {};
				\node[gvertex] (5) at (0.27561, 0.70119) {};
				\draw[edge] (0) to (1) (0) to (2) (0) to (3) (1) to (2) (1) to (4) (2) to (5) (3) to (4) (3) to (5) (4) to (5);
			\end{tikzpicture}%
		};
		\node[elem,font=\tiny,anchor=south west] (code) at (0,-11.75) {%
			\begin{minipage}{6.2cm}
				\begin{pythonfig}[tabsize=2,gobble=10]
					def is_min_rigid(
						self,
						dim: int = 2,
						algorithm: str = "default",
						...
					) -> bool:
						"""
						Return whether the graph is minimally ``dim``-rigid.

						Definitions
						-----------
						:prf:ref:`Minimal dim-rigidity <def-min-rigid-graph>`

						Parameters
						----------
						dim:
							Dimension.
						algorithm:
							...
							If ``"sparsity"`` (only if ``dim=2``),
							then :prf:ref:`(2,3)-tightness
							<def-kl-sparse-tight>` and
							:prf:ref:`thm-2-gen-rigidity` are used.
							...
				\end{pythonfig}\vspace{-0.25cm}%
			\end{minipage}%
		};
		\node[elem,font=\scriptsize,anchor=south east,align=left] (out) at (13.5,-11.75) {%
			\begin{minipage}{4cm}
				\begin{pythonfig}[tabsize=2,gobble=10]
					from pyrigi import graphDB
					G = graphDB.ThreePrism()
					G.is_min_rigid(dim=2)
				\end{pythonfig}\vspace{-0.25cm}%
			\end{minipage}%
		};
		\node[head] at ($(thm.north west)+(0.25,0)$) {Math Documentation};
		\node[head] at ($(doc.north west)+(0.25,0)$) {Math Documentation};
		\node[head] at ($(cdoc.north west)+(0.25,0)$) {Code Documentation};
		\node[head] at ($(ex.north west)+(0.25,0)$) {Graph};
		\node[head] at ($(code.north west)+(0.25,0)$) {Code};
		\node[head] at ($(out.north west)+(0.25,0)$) {Usage};
		\node[labelsty,opacity=0.4] (fdoc) at (6.75,0) {further doc};

		\draw[ref] ($(code.south east)+(0,0.75)+(-1.5,0)$) -| ++(2,0) |- ($(thm.north west)+(0,-0.5)$);
		\draw[ref] ($(thm.north west)+(0,-1.67)+(0.5,0)$) -- ++(-2,0) |- ($(doc.north east)+(0,-0.5)$);
		\draw[ref] ($(code.north west)+(0,-3.39)+(0.4,0)$) -- ++(-0.85,0) |- ($(doc.north west)+(0,-0.5)$);
		\draw[ref] ($(doc.south east)+(0,0.875)+(-2.25,0)$) -- ++(3,0) |- ($(cdoc.north west)+(0,-0.5)$);
		\draw[sref] ($(doc.south east)+(0,1.5)+(-0.5,0)$) -| ($(fdoc.south)+(0,0)$);
		\draw[sref] ($(doc.south east)+(0,1.825)+(-1,0)$) -| ($(fdoc.south)+(-0.25,0)$);
		\draw[srefrev] ($(fdoc.south)+(0.25,0)$) |-  ($(thm.south west)+(0,1.475)+(0.5,0)$);
		\draw[srefrev] ($(fdoc.south)+(0.25,0)$) |- ($(code.south east)+(0,1.27)+(-2.417,0)$);
		\draw[sref] ($(thm.south east)+(0,0.88)+(-2,0)$) -| ++(2.4,0.5) node[above,anchor=west,black,labelsty,opacity=0.4,rotate=90] (ref) {bibliography};

		\draw[con] ($(out.north east)+(0,-0.9)$) -- ++(0.5,0) |- node[pos=0.25,above,anchor=south,labelsty,rotate=90] {takes example} (ex.east);
		\draw[con] ($(code.north west)+(3,0)$) to node[right,labelsty] {implements} ($(doc.south west)+(3,0)$);
		\draw[con] (code.east) -- ++(1.25,0) |-node[pos=0.1,below,labelsty,rotate=90] {produces} (cdoc.west);
		\draw[conb] ($(out.south west)+(0,0.6)$) to ++(-1.25,0) |- ($(code.north east)-(0,0.5)$);
		\draw[con] ($(out.south west)+(0,0.6)$) to node[below,labelsty] {uses} ++(-1.25,0) |- ($(code.north east)-(0,0.5)$);
	\end{tikzpicture}
	\caption{Example structure of graph rigidity in \pyrigi.
		The rounded rectangles describe different parts of \pyrigi.
		Blue dashed arrows \protect\tikz{\protect\draw[ref](0,0)--(0.5,0);} symbolize references between elements
		while red arrows \protect\tikz{\protect\draw[con](0,0)--(0.5,0);} show structural relations.}
	\label{fig:pyrigi-graph}
\end{figure}

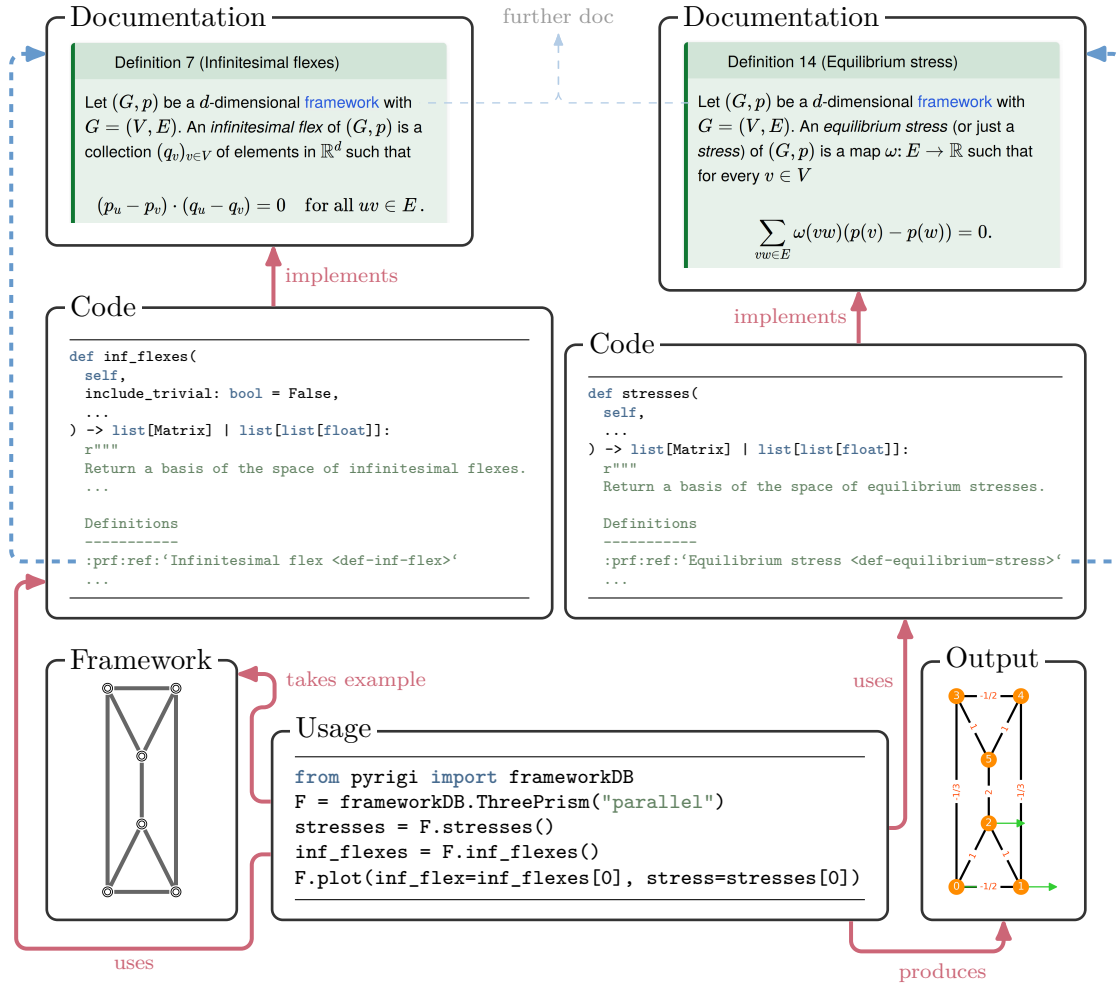
\begin{figure}[ht]
	\centering
	\begin{tikzpicture}
		\node[elem,anchor=north west] (doc1) at (0,0) {\includegraphics[width=5cm]{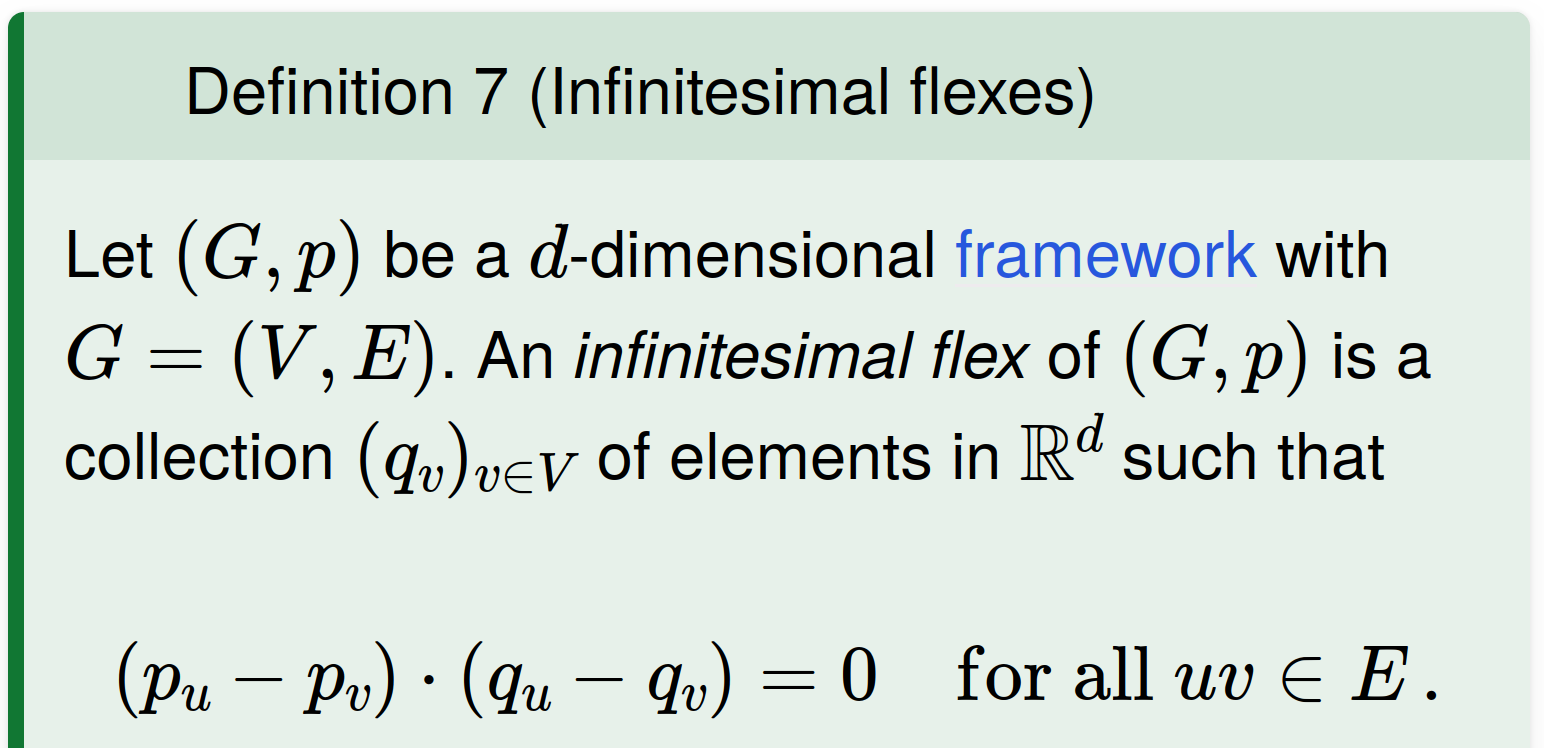}};
		\node[elem,anchor=north east] (doc2) at (13.7,0) {\includegraphics[width=5cm]{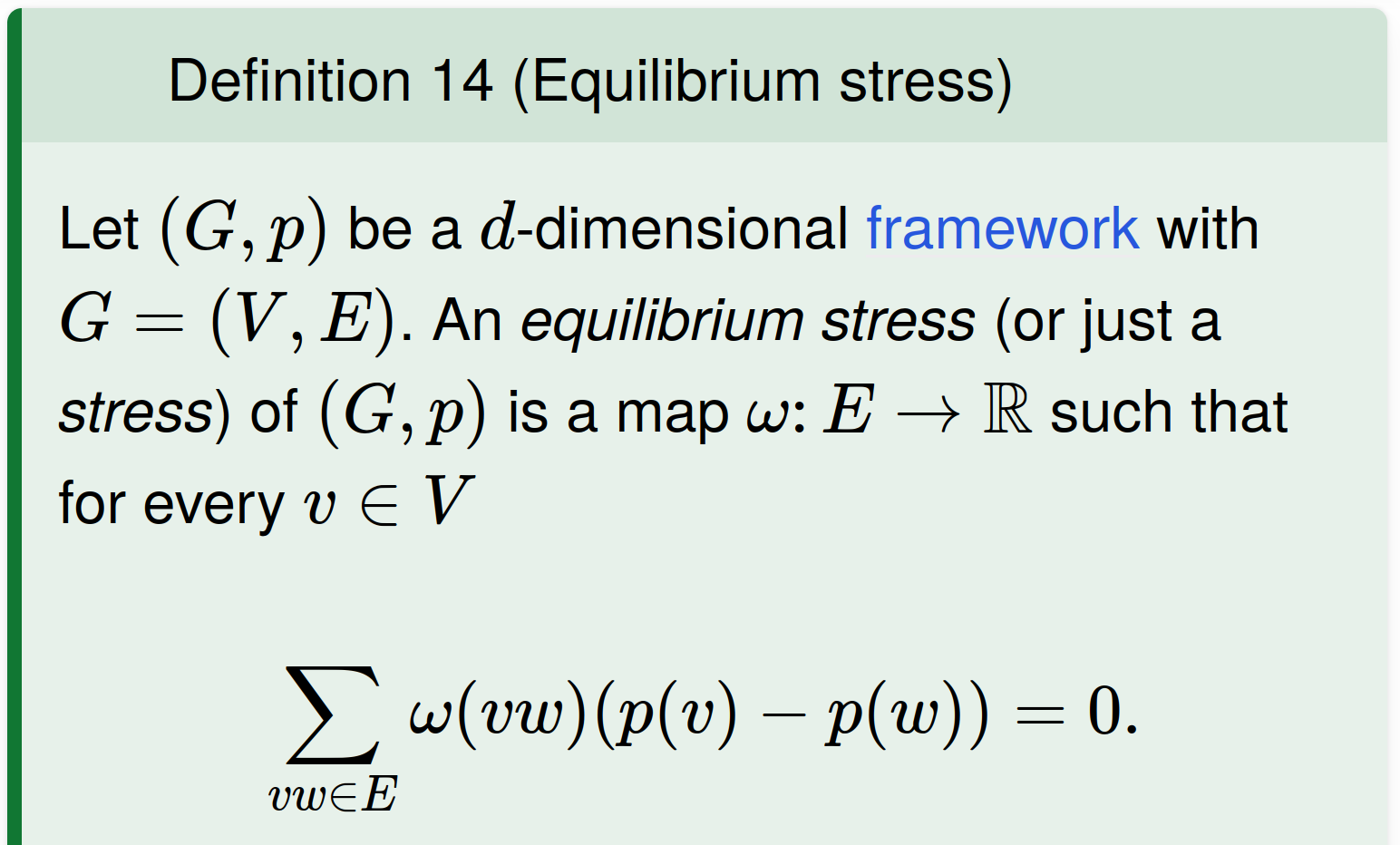}};
		\node[elem,anchor=south west,minimum size=2.5cm] (ex) at (0,-12) {%
			\begin{tikzpicture}[%
				scale=0.45,
				fvertex/.style={circle,inner sep=0pt,minimum size=3pt,fill=white,draw=black,double=white,double distance=0.25pt,outer sep=1pt},
				edge/.style={line width=1.5pt,black!60!white,}
			]
				\node[fvertex] (0) at (0, 0) {};
				\node[fvertex] (1) at (2, 0) {};
				\node[fvertex] (2) at (1, 2) {};
				\node[fvertex] (3) at (0, 6) {};
				\node[fvertex] (4) at (2, 6) {};
				\node[fvertex] (5) at (1, 4) {};
				\draw[edge] (0) to (1) (0) to (2) (0) to (3) (1) to (2) (1) to (4) (2) to (5) (3) to (4) (3) to (5) (4) to (5);
			\end{tikzpicture}%
		};
		\node[elem,font=\tiny,anchor=south west] (code1) at (0,-8) {%
			\begin{minipage}{6.05cm}
				\begin{pythonfig}[tabsize=2,gobble=10]
					def inf_flexes(
						self,
						include_trivial: bool = False,
						...
					) -> list[Matrix] | list[list[float]]:
						r"""
						Return a basis of the space of infinitesimal flexes.
						...

						Definitions
						-----------
						:prf:ref:`Infinitesimal flex <def-inf-flex>`
						...
				\end{pythonfig}\vspace{-0.25cm}%
			\end{minipage}%
		};
		\node[elem,font=\tiny,anchor=south east] (code2) at (13.7,-8) {%
			\begin{minipage}{6.24cm}
				\begin{pythonfig}[tabsize=2,gobble=10]
					def stresses(
						self,
						...
					) -> list[Matrix] | list[list[float]]:
						r"""
						Return a basis of the space of equilibrium stresses.

						Definitions
						-----------
						:prf:ref:`Equilibrium stress <def-equilibrium-stress>`
						...
				\end{pythonfig}\vspace{-0.25cm}%
			\end{minipage}%
		};
		\node[elem,font=\scriptsize,anchor=south east,align=left] (out) at (11.1,-12) {%
			\begin{minipage}{7.5cm}
				\begin{pythonfig}[tabsize=2,gobble=10]
					from pyrigi import frameworkDB
					F = frameworkDB.ThreePrism("parallel")
					stresses = F.stresses()
					inf_flexes = F.inf_flexes()
					F.plot(inf_flex=inf_flexes[0], stress=stresses[0])
				\end{pythonfig}\vspace{-0.25cm}%
			\end{minipage}%
		};
		\node[elem,anchor=south east] (outf) at (13.7,-12) {\includegraphics[width=1.55cm]{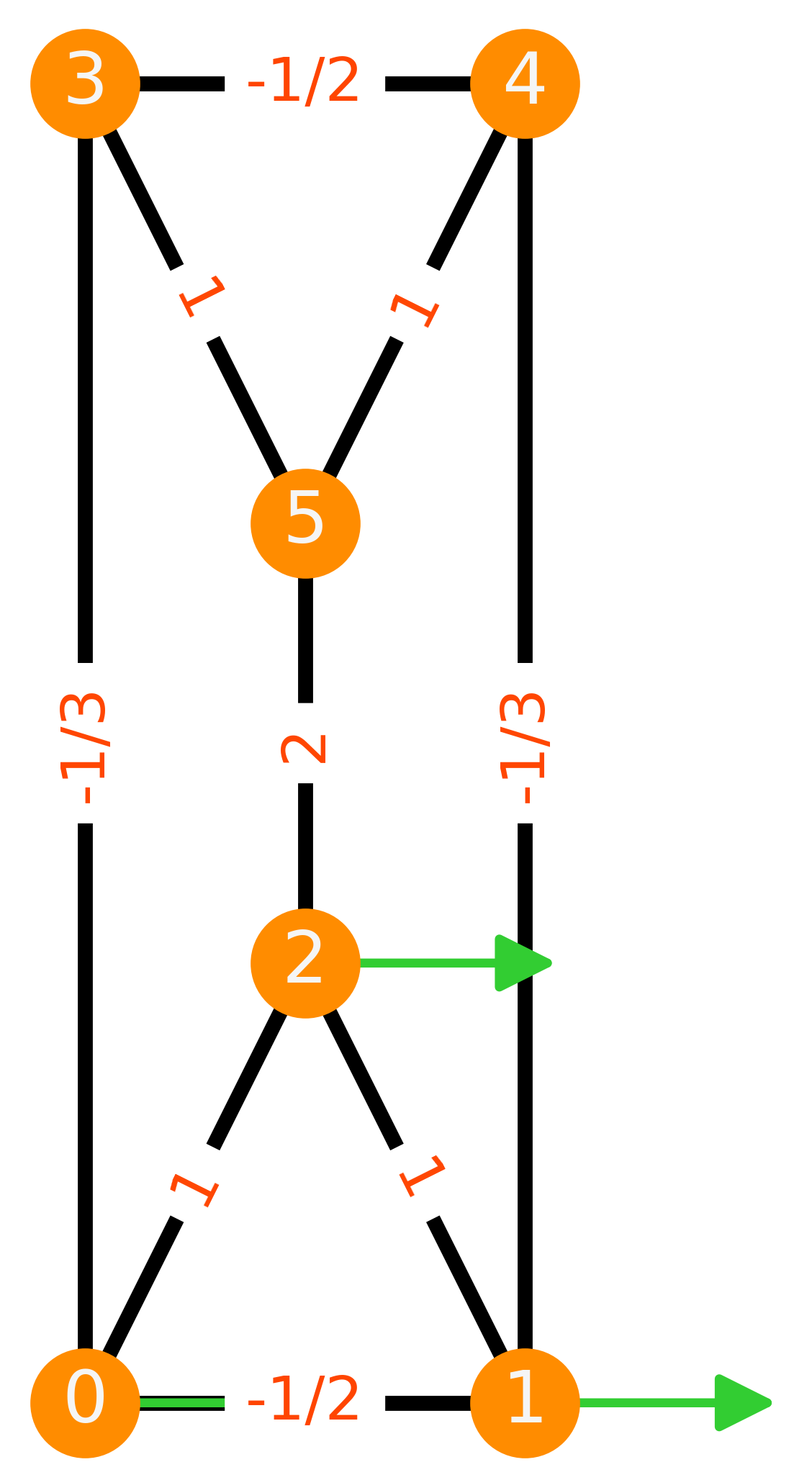}};
		\node[head] at ($(doc1.north west)+(0.25,0)$) {Documentation};
		\node[head] at ($(doc2.north west)+(0.25,0)$) {Documentation};
		\node[head] at ($(ex.north west)+(0.25,0)$) {Framework};
		\node[head] at ($(code1.north west)+(0.25,0)$) {Code};
		\node[head] at ($(code2.north west)+(0.25,0)$) {Code};
		\node[head] at ($(out.north west)+(0.25,0)$) {Usage};
		\node[head] at ($(outf.north west)+(0.25,0)$) {Output};
		\node[labelsty,opacity=0.4] (fdoc) at (6.75,0) {further doc};

		\draw[ref] ($(code1.south west)+(0,0.77)+(0.4,0)$) -- ++(-0.85,0) |- ($(doc1.north west)+(0,-0.5)$);
		\draw[ref] ($(code2.south east)+(0,0.77)+(-0.2,0)$) -- ++(0.675,0) |- ($(doc2.north east)+(0,-0.5)$);
		\draw[srefrev] ($(fdoc.south)+(0,0)$) |- ($(doc1.north east)+(0,-1.15)+(-0.71,0)$);
		\draw[srefrev] ($(fdoc.south)+(0,0)$) |- ($(doc2.north west)+(0,-1.15)+(0.48,0)$);

		\draw[con] ($(out.north east)+(0,-1.3)$) -| node[pos=0.85,left,labelsty] {uses} ($(code2.south)+(1.05,0)$);
		\draw[con] ($(out.north west)+(0,-0.95)$) -- ++(-0.25,0) -- ++(0,1.25) --++(0.3,0) |- node[pos=0.4,right,anchor=west,labelsty,align=center] {takes example} ($(ex.north east)+(0,-0.2)$);
		\draw[con] ($(code1.north west)+(3,0)$) to node[right,labelsty] {implements} ($(doc1.south west)+(3,0)$);
		\draw[con] ($(code2.north east)+(-3,0)$) to node[left,labelsty] {implements} ($(doc2.south east)+(-3,0)$);
		\draw[con] ($(out.north west)+(0,-1.65)$) -- ++(-0.25,0) -- ++(0,-1.25) to node[below,labelsty] {uses} ++(-3.1,0) |-  ($(code1.south west)+(0,0.5)$);
		\draw[con] ($(out.south east)+(-0.5,0)$) -- ++(0,-0.4) -| node[pos=0.3,below,labelsty] {produces} (outf.south);
	\end{tikzpicture}
	\caption{Framework rigidity in \pyrigi.
		The rounded rectangles describe different parts of \pyrigi. Arrow styles are as in \Cref{fig:pyrigi-graph}.}
\label{fig:pyrigi-framework}
\end{figure}

A basic part of a software documentation is the description of the code and the functions and methods that are available.
\pyrigi\ goes a step further, adding a second part that combines the mathematical concepts and theorems of rigidity theory with computational tools.
This means that each method and function in \pyrigi\ --- unless it is already well-known --- is linked to its respective mathematical definition. For algorithmic concepts and implementation decisions there is usually a mathematical theorem and
an algorithm description in the documentation.
The mathematical documentation collects all necessary information in one place in a handbook-like style organized in topical pages.
Definitions and theorems can therefore be referenced both from the code and the documentation.
This provides the user with the necessary knowledge about what is computed and how this is done.
The interactions between the theoretical and the computational side are visualized in \Cref{fig:pyrigi-graph,fig:pyrigi-framework} for a graph and a framework example, respectively.
In the following paragraphs, we describe these interactions in more detail.
Besides the in-depth documentation which is intended for more advanced users, there are of course also tutorials for getting started.

Let us start by discussing the interplay of mathematical concepts, documentation, and code on an exemplary graph depicted in \Cref{fig:pyrigi-graph}.
A crucial definition in rigidity theory is that of a minimally rigid graph.
Out of the many possible equivalent definitions, we chose to adopt one related to the rigidity of frameworks (top left of the figure).
However, the genericity that is mentioned in this definition is, from a coding point of view, difficult to use directly.
The theorem of Pollaczek-Geiringer and Laman (top right of the figure) yields a better strategy for computation.
Bibliographic references for the theorem guarantee scientific standards of the mathematical documentation.
The theorem links to another concept of graph theory, namely tightness, for which there is (for some important cases) a polynomial time algorithm.
This algorithm is one option to check rigidity in the method \pythoninline{is_min_rigid} as can be seen in the docstring of the method (bottom left of the figure).
The docstring also refers to the respective definitions and theorems.
Additionally, the docstring shows an explanation of the method and the parameters and it gives examples for their usage, fulfilling common coding standards, such as PEP 8.
All of this is also presented in the online documentation in a well-formatted and accessible manner.
Note that in the figure some parts are left out (indicated by ``\ldots'') for space and illustration purposes.
The example in the figure is taking a predefined graph from the database of graphs \pythoninline{graphDB} that accompanies \pyrigi.

The example structure for framework (\Cref{fig:pyrigi-framework}) shows two definitions relevant for frameworks, namely infinitesimal flexes and stresses (top of the figure).
Frameworks are implemented in \pyrigi\ as their own class.
Two of its methods are shown: \pythoninline{inf_flexes} computes all infinitesimal flexes of a given framework and \pythoninline{stresses} does so with the equilibrium stresses (center of the figure).
References to definitions are considered an important aspect of the documentation in \pyrigi\ and hence there is a standard section for them in the docstrings, as can be seen in the figure.
The example is chosen to show that computational results of frameworks can easily be visualized in \pyrigi.
The input framework is obtained from the framework database \pythoninline{frameworkDB} that is part of \pyrigi.
The output figure shows both the infinitesimal flexes (as green arrows) as well as the stresses (as labels on the edges).

Besides the documentation on the code and the mathematical background, \pyrigi{}'s documentation also provides a tutorial on the rigidity of graphs and frameworks\footnote{\url{https://pyrigi.github.io/PyRigi/userguide/tutorials/rigidity.html}}.
This tutorial serves as an introduction to the functionality of \pyrigi.
It contains a description and the usage of some basic commands, including the elements described in the figures.
Again, references to the mathematical definitions and the method descriptions help the user to easily navigate through the documentation.

\section{Functionality}
\label{functionality}

In this section, we illustrate the current functionality of \pyrigi.
We emphasize what \pyrigi\ can do and refer the reader to the online documentation of \pyrigi\ to see how this can be used in practice.
We start by describing how frameworks are handled by \pyrigi{},
and then move to graphs, concluding with motions.
We chose this order, since many notions related to graph rigidity
are first introduced for frameworks and later defined for graphs
by considering so-called ``generic'' realizations.

\subsection{Frameworks}
\label{functionality_framework}
Let $G=(V,E)$ be a simple graph, $d\in \mathbb{N}$ and let $p:V\rightarrow \mathbb{R}^d$
be a \emph{$d$-dimensional realization of $G$}.
A \emph{$d$-dimensional framework} is then defined as the pair $(G,p)$.
In \pyrigi{}, a framework object can be instantiated using the class \pythoninline{Framework}.
It expects a \pythoninline{Graph} and a realization given by a dictionary
as input. For instance, a framework on a 4-cycle can be constructed by calling

\begin{python}
from pyrigi import Framework, Graph
G = Graph([(0,1), (1,2), (2,3), (3,0)])
realization = {0: [0,0], 1: ["sqrt(2)",0], 2: [1,1], 3: [0,"3/4"]}
F = Framework(G, realization)
\end{python}
The coordinates in the realization can be provided as \pythoninline{sympy} expressions using their string representation,
so that the computations on the framework can be performed symbolically.
Specific relevant frameworks are stored in the internal \pythoninline{frameworkDB}.
For instance, we can construct a parallel configuration of the 3-prism using the following commands:

\begin{python}
from pyrigi import frameworkDB
TP = frameworkDB.ThreePrism(realization="parallel")
\end{python}
This framework consists of two flipped and translated triangles,
with the property that corresponding vertices are connected by edges, that lie on three parallel lines. It is depicted in \Cref{fig:pyrigi-framework}.
We use the frameworks \pythoninline{F} and \pythoninline{TP} as examples throughout this subsection.

\minisec{Equivalence and Congruence}
To distinguish whether a given framework is flexible or rigid and to say that one framework is a deformed version of another framework, we employ a theoretical notion to compare realizations of the same graph up to rigid motions.
Two $d$-dimensional frameworks $(G,p)$ and $(G,p')$ on a graph $G=(V,E)$
are called \emph{equivalent} if
\begin{align*}
	||p(u)-p(v)|| &= ||p'(u)-p'(v)|| &&\text{for all } uv\in E.\\
	\intertext{They are called \emph{congruent} if}
	||p(u)-p(v)|| &= ||p'(u)-p'(v)|| &&\text{for all } u,v\in V.
\end{align*}
In \pyrigi, these concepts can be checked by calling the methods \pythoninline{is_equivalent}
and \pythoninline{is_congruent}.
By default, these checks are performed symbolically.
When faster computation is required, users may switch to numerical methods by setting the appropriate keyword.
This converts the internal coordinates to floating point numbers and runs numerical checks instead of symbolic ones.
The computation uses a predefined numeric tolerance, which also can be specified.

We then say that a $d$-dimensional framework $(G,p)$ is \emph{(continuously) rigid} if every framework $(G, q)$ equivalent to $(G, p)$, and sufficiently close to it in the Euclidean distance, is actually congruent to it; otherwise, it is called \emph{(continuously) flexible}.

\minisec{Transformations}

As we are often interested in comparing properties of different frameworks, \pyrigi\ comes with various tools for transforming frameworks.
First and foremost, we can modify the underlying graph by adding and deleting vertices or edges.
We can project a realization to a lower dimension by calling \pythoninline{projected_realization}.
In addition, it is possible to \pythoninline{rescale} all coordinates by a factor, we can \pythoninline{translate} a framework and it is possible to \pythoninline{rotate} it in $\mathbb{R}^2$ and $\mathbb{R}^3$.

\minisec{Infinitesimal Rigidity}
Determining whether a given framework with arbitrary coordinates is rigid in $\mathbb{R}^d$ is coNP-hard for $d\geq 2$ \parencite{Abbott2008}, meaning that checking whether a framework is flexible is NP-hard.
Therefore, stronger conditions are typically employed to determine whether a framework is rigid.

The concept of infinitesimal rigidity provides an important sufficient condition that arises from showing that a framework's realization is a non-singular real isolated solution of the underlying Euclidean distance constraint system. In other words, a $d$-dimensional framework $(G,p)$ on a graph $G=(V,E)$ is \emph{infinitesimally rigid}, if the rank of the Jacobian corresponding to the polynomial system of the associated edge-length equations --- also referred to as the \emph{rigidity matrix} --- is equal to $d|V|-\binom{d+1}{2}$. This matrix is the $|E| \times d|V|$ matrix whose row labelled by $e = \{u,v\}$ is
\[
	\begin{array}{ccccccccccc}
		&&& u &&&& v \\
		\bigl( 0 & \cdots & 0 & p(u) - p(v) & 0 & \cdots & 0 & p(v) - p(u) & 0 & \cdots & 0 \bigr) \,.
	\end{array}
\]

It turns out that the previously constructed framework \pythoninline{F} on the 4-cycle graph is not infinitesimally rigid. When adding the edge $\{1,3\}$ by calling \pythoninline{F.add_edge([1,3])}, it becomes a Diamond framework which is \emph{minimally infinitesimally rigid}, meaning that the removal of any edge from the framework's underlying graph causes it to become infinitesimally flexible.
This fact can be verified using the command \pythoninline{F.is_min_inf_rigid()}.
If we further add the edge~$\{0,2\}$,
the framework \pythoninline{F} becomes a realization of the complete graph on 4 vertices $K_4$.
It is \emph{redundantly infinitesimally rigid}, meaning that the removal of any edge still produces an infinitesimally rigid framework,
as can be verified by calling \pythoninline{F.is_redundantly_inf_rigid()}.
The keyword \pythoninline{numerical} can be used in all of these cases to speed up the computations
by converting the symbolic coordinates to floating point numbers and computing the rigidity matrix's kernel using \emph{NumPy}~\parencite{numpy}.

\minisec{Stresses and Flexes}
Conversely, the framework \pythoninline{TP} on the 3-prism graph is not infinitesimally rigid,
as can be checked via \pythoninline{TP.is_inf_rigid()}, which returns
\pythoninline{False}.

We can verify this result by computing the framework's \emph{infinitesimal flexes},
which are given by the kernel of the rigidity matrix. In other words,
the infinitesimal flexes of framework $(G,p)$ are vectors $(q_v)_{v\in V}$
at each vertex satisfying $(p_u-p_v)\cdot (q_u-q_v)=0$ for each edge $\{u,v\}\in E$ and $p_v=p(v)$.
They are called \emph{trivial}, if they extend to an ambient isometry, \emph{nontrivial} otherwise.
In particular, an infinitesimally rigid framework possesses no nontrivial infinitesimal flexes.
For instance, all infinitesimal flexes of the 3-prism framework \pythoninline{TP} can be computed by calling \pythoninline{TP.inf_flexes(include_trivial=True)}.

An associated concept are \emph{equilibrium stresses}, which form the cokernel of the rigidity matrix.
Computing \pythoninline{inf_flexes=TP.inf_flexes()} and \pythoninline{stresses=TP.stresses()} demonstrates that the framework has both, a one-dimensional space of nontrivial infinitesimal flexes and equilibrium stresses.
We can then verify that these flexes and stresses are indeed nontrivial infinitesimal flexes and equilibrium
stresses by calling \pythoninline{TP.is_nontrivial_flex(inf_flexes[0])} and \pythoninline{TP.is_stress(stresses[0])},
respectively. These computations can be performed numerically as well.

\minisec{Prestress Stability and Second-Order Rigidity}
Even though the parallel realization of the 3-prism graph is not infinitesimally rigid,
the framework is prestress stable.
This concept constitutes a sufficient criterion, weaker than infinitesimal rigidity,
for the continuous rigidity of a framework.
It depends on the idea that stresses may \emph{block} certain infinitesimal flexes from extending to continuous motions.
More precisely, if there exists an equilibrium stress $\omega$ such that for every nontrivial infinitesimal motion $q$ it holds that
\[\sum_{ij\in E} \omega_{ij} \cdot ||q_i-q_j||^2 > 0,\]
the framework is called \emph{prestress stable} \parencite{ConnellyWhiteley1996}.
For the 3-prism framework \pythoninline{TP}, this can be checked by calling \pythoninline{TP.is_prestress_stable()}
which returns \pythoninline{True}, so it is indeed rigid.
If the quantifiers are reversed, we obtain the weaker concept of \emph{second-order rigidity} \parencite{ConnellyWhiteley1996}.
In other words, this is the case if for every nontrivial infinitesimal motion $q$ there exists an equilibrium stress $\omega$ such that the inequality above holds.
This criterion can be verified via the method \pythoninline{is_second_order_rigid}.
Since prestress stability and second-order rigidity both imply (continuous) rigidity \parencite{ConnellyGortler2017}, the framework cannot be continuously deformed.
This check is based on the sums of nonnegative circuits decomposition \parencite{IlimandeWolff2016}, which is a general sufficient condition for the nonnegativitiy (or positivity) of polynomials.
It is currently required that either the space of infinitesimal flexes or the space of equilibrium stresses is 1-dimensional, though we intend to extend this approach to more general systems in the near future.
It is possible to use the \pythoninline{numerical} keyword here as well.

\minisec{Drawing and Plotting}
In \pyrigi{}, it is further possible to input a graph via the provided \pythoninline{GraphDrawer} class and simultaneously provide it with a realization.
The GUI can be opened, for instance, in a Jupyter Notebook with the command

\begin{figure}[ht]
	\centering
	\begin{minipage}[c]{0.66\linewidth}
		\includegraphics[width=1\linewidth]{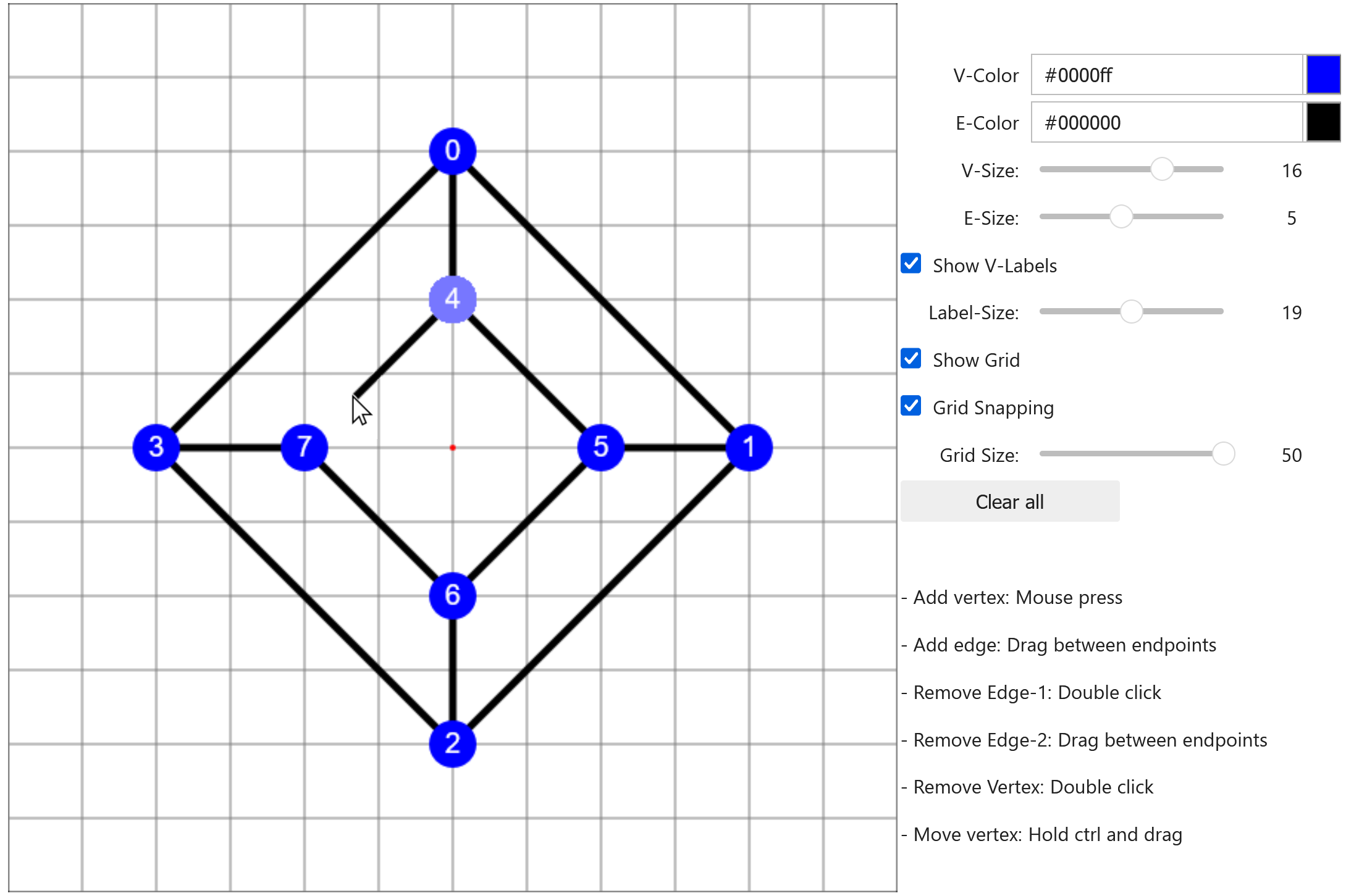}
	\end{minipage}\hfill
	\begin{minipage}[c]{0.32\linewidth}
		\includegraphics[width=1\linewidth]{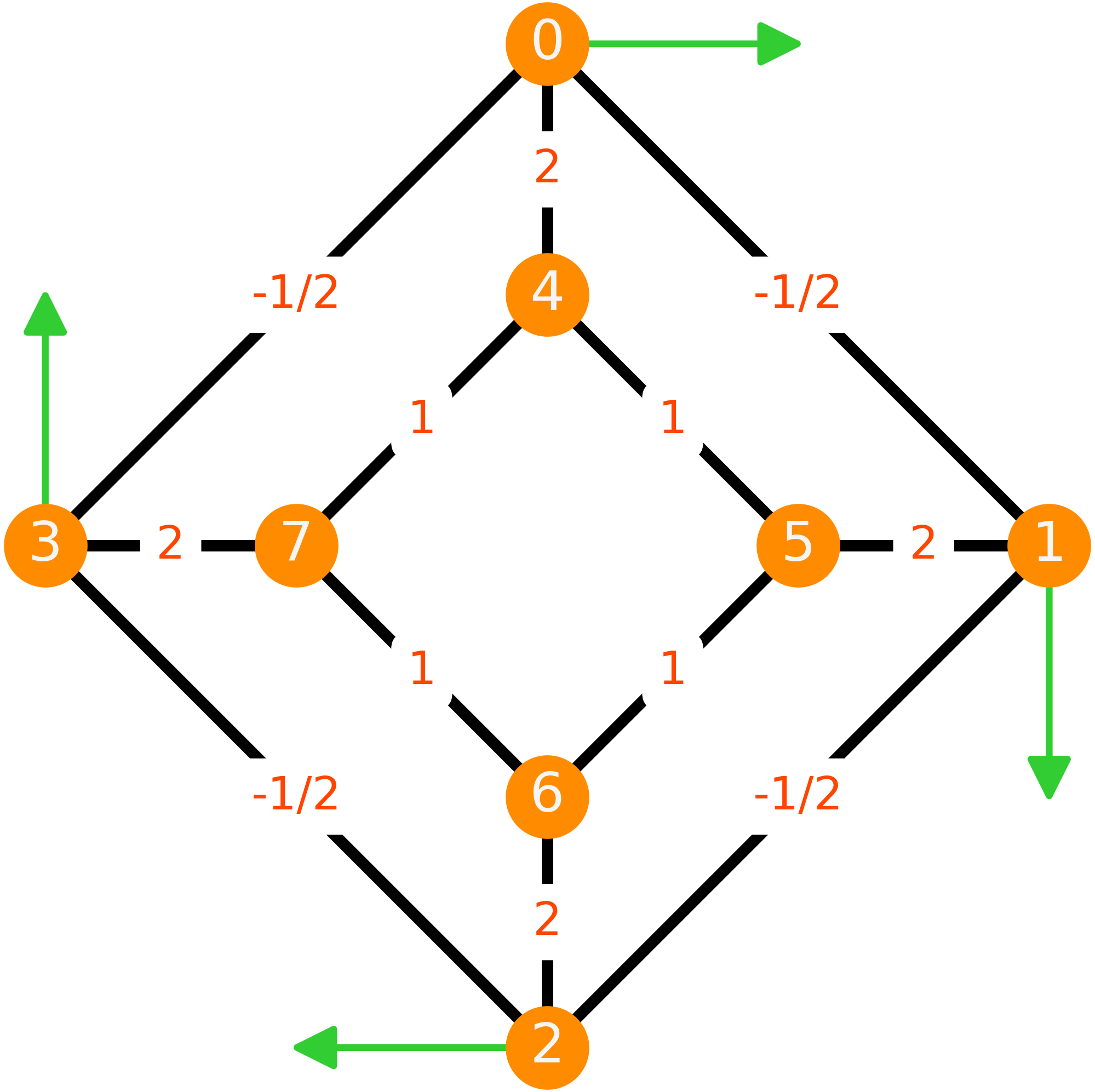}
	\end{minipage}
	\caption{The \texttt{GraphDrawer} GUI makes it possible to input a graph with a realization in \pyrigi{} by drawing it (left).
		In this screenshot from a Jupyter Notebook, the edge $\{4,7\}$ is currently drawn.
		The resulting framework is then visualized using the internal plotting routines,
		while showing an infinitesimal flex and an equilibrium stress (right).}
	\label{fig:graphdrawer}
\end{figure}

\begin{python}
from pyrigi import GraphDrawer
GD = GraphDrawer()
\end{python}
and an exemplary canvas is depicted in \Cref{fig:graphdrawer}~(left).

Interaction with the GUI is driven by left-clicking the canvas to create vertices that are numbered in ascending order.
Furthermore, we can create edges by clicking and holding an existing vertex and dragging the cursor to another vertex. Conversely, vertices and edges can be removed by double-clicking.
Additionally, it is possible to add a grid to the canvas of a chosen fineness and to snap the vertices to the grid.
This allows one for a better control of the vertex positions.
The created graph and framework can be cast to \pyrigi{} objects by calling methods \pythoninline{graph} and \pythoninline{framework}, respectively.
For example, the output from \Cref{fig:graphdrawer}~(left) can be obtained with the command \pythoninline{F4 = GD.framework()}.

Moreover, \pyrigi's \pythoninline{Framework} class comes with a variety of visualization options.
By calling \pythoninline{plot2D} or \pythoninline{plot3D}, a given framework can be depicted
in $\mathbb{R}^2$ or $\mathbb{R}^3$, respectively.
Frameworks in $d$-dimensional space with large $d\in \mathbb{N}$
can be projected to 2D or 3D by providing a projection matrix with appropriate dimensions.
Additionally, it is possible to depict infinitesimal flexes and equilibrium stresses,
see \Cref{fig:graphdrawer}~(right).
Regarding style adjustments, for instance edges can be depicted as curved arcs
or can be assigned custom edge colors.
The edge colors can also be provided as a partition of the underlying graph's edges
which automatically creates distinct colors for each partition suitable for red–green color blind people.
These parameters can be saved in an object of the \pythoninline{PlotStyle}
class so that the particular style becomes reproducible across multiple instances of frameworks.
Many additional keywords from the \emph{Matplotlib} library can be used here, as well.

\minisec{Exports}

\pyrigi{} makes it possible to export frameworks in two separate formats.
First, the method \pythoninline{to_tikz} creates the TikZ\footnote{\url{https://github.com/pgf-tikz/pgf}} code to embed the framework in \pyrigi{}'s style in a scientific article that is written in \LaTeX. Second, it is possible to create an STL file that can be used to 3D print the framework via the method \pythoninline{generate_stl_bars}, which renders simple bars with holes that can be connected by screws and stop-nuts.

\subsection{Graphs}

Graphs are one of the main objects in \pyrigi.
A graph object can be constructed using the class \pythoninline{Graph}.
The most common way of constructing a graph is by providing a list of edges of the graph
as we have seen at the beginning of \Cref{functionality_framework}.
There are also other ways of constructing a \pythoninline{Graph},
as for example via the method \pythoninline{from_vertices_and_edges} that accepts a pair of lists indicating the vertices and the edges.
This is indicated if one wants to have isolated vertices.

\minisec{Generic Rigidity}

As one may expect, \pyrigi{} can handle generic $d$-rigidity of graphs,
where $d$ denotes the dimension,
with algorithms that are tailored to the considered dimension.
Here, by \emph{generic $d$-rigidity} we mean the following.
A graph $G$ is called \emph{$d$-rigid} if any $d$-dimensional framework $(G, p)$
such that all entries of the image of $p$ are algebraically independent
(i.e., the realization $p$ is \emph{generic}) is infinitesimally rigid.
Otherwise, the graph is \emph{$d$-flexible}.

Here, the most important method is \pythoninline{is_rigid}.
In dimension~$2$, \pyrigi{} adopts by default an algorithm,
that checks whether there exists a $(2,3)$-tight subgraph of the given graph with the same vertex set; see \textcite{PollaczekGeiringer1927, Laman1970,LeeStreinu2008}.
Here, a $G = (V, E)$ is \emph{$(2,3)$-tight} if $|E| = 2|V| - 3$
and it is \emph{$(2,3)$-sparse}, that means
every subset of $n'$ vertices of $G$ that determine at least one edge, spans at most $2n' - 3$ edges in~$G$.

In dimension $1$, instead, \pyrigi{} uses by default the well-known characterization
stating that a graph is rigid if and only if it is connected.
In addition to these two algorithms, for every dimension \pyrigi{} may adopt
a randomized method (using parameter \pythoninline{algorithm="randomized"}),
based on creating a framework with the given underlying graph
whose realization has randomly chosen integer coordinates,
and checking whether the rigidity matrix kernel has the expected dimension;
see \textcite[Alg~5.2, Prop~5.7]{GortlerHealyThurston2010}.
This randomized method never gives false positives.
The theoretical foundations of this algorithm give a precise estimate of the probability of obtaining a false negative
in terms of the range in which the randomly coordinates are chosen.
By this it is possible to specify in \pyrigi\ an upper bound for this probability when calling the algorithm.

Similarly to \pythoninline{is\_rigid},
one may use the method \pythoninline{is\_min\_rigid}
to check whether a graph is minimally rigid.
As a step further, \pyrigi{} is able to determine the rigid components of a graph,
which are the maximal rigid subgraphs:
the method \pythoninline{rigid\_components} returns a list of their vertices.
Finally, the number of complex realizations in the plane of a minimally rigid graph
up to the action of the complexification of the group of isometries of~$\mathbb{R}^2$
can be computed via the method \pythoninline{number_of_realizations}
which is a wrapper for the analogous function from the package \pythoninline{lnumber} \parencite{Capco2024};
see \textcite{CapcoEtAl2018} for the theoretical background.
If the parameter \pythoninline{spherical} is set to \pythoninline{True},
the algorithm computes the number of complex realizations on the sphere.

\minisec{$(k,\ell)$-Sparsity and Tightness}

As already mentioned while discussing generic rigidity, the notions of \emph{$(k, \ell)$-sparsity} and \emph{tightness} play a prominent role in rigidity theory,
and \pyrigi{} ships efficient methods for checking both of them,
based on the so-called \emph{pebble game} algorithm \parencite{JacobsHendrickson1997,LeeStreinu2008}.
Here, $(k, \ell)$-sparsity and tightness are generalizations of the $(2,3)$-concepts introduced before
and are defined analogously.

\minisec{Redundant Rigidity}

We say that a graph $G$ is \emph{redundantly $d$-rigid} if removing any edge from $G$ yields a $d$-rigid graph.
In a similar manner, one can define \emph{vertex redundantly $d$-rigidity} and the counterparts of both these concepts involving removing $k$ edges or vertices.
\pyrigi{} has implemented methods for checking all these properties,
based on the works of \textcite{Servatius1989,YuAnderson2009,SummersYuAnderson2009,KaszanitzkyKiraly2016,Jordan2016,Jordan2021,JordanPostonRoach2022,MotevallianYuAnderson2015}.

\minisec{Global Rigidity}

A framework $(G, p)$ is called \emph{globally $d$-rigid}
if every $d$-dimensional realization $q$ that is equivalent to $p$ is also congruent to $p$.

A celebrated theorem by \textcite{GortlerHealyThurston2010}
shows that either all frameworks $(G, p)$, when $p$ is generic, are globally $d$-rigid, or none is so.
Therefore, in this way one can speak of global rigidity as a property of a graph.
In the same paper, \citeauthor{GortlerHealyThurston2010} provide a probabilistic algorithm to check global rigidity of graphs, which is implemented in \pyrigi{}.
Detecting whether a graph is globally rigid or not is implemented in \pyrigi{} via the method \pythoninline{is\_globally\_rigid}.
As for the case of \pythoninline{is_rigid}, the user may enter a parameter
serving as an upper bound for the probability of false negatives.
In dimension~$2$, \pyrigi{} relies also on the result by \textcite{JacksonJordan2005} stating that
a graph is globally $2$-rigid if and only it is $3$-connected and redundantly rigid;
a user may access this functionality by calling the method \pythoninline{is\_globally\_rigid}
with the parameter \pythoninline{algorithm="redundancy"}.

A concept related to global rigidity is the one of \emph{global linkedness}.
If $(G, p)$ is a framework and $u,v$ are vertices of $G$, we say that $u$ and $v$ are globally linked in $(G, p)$
if the distance between $p(u)$ and $p(v)$ equals the distance between $q(u)$ and $q(v)$
for any other realization $q$ of $G$ equivalent to $p$.
Then, two vertices $u$ and $v$ of a graph $G$ are said to be \emph{weakly globally linked}
if they are globally linked for at least one generic realization of $G$.
This notion turned out to be useful in proving the generalization of the celebrated result of Lov\`asz and Yemini
stating that every $d(d+1)$-connected graph is globally rigid in $\mathbb{R}^d$ \parencite{Villanyi2025}.
\textcite{JordanVillanyi2024} provide an algorithm to check whether a pair of vertices is weakly globally linked.
This algorithm is implemented in \pyrigi{} in the method \pythoninline{is_weakly_globally_rigid}.

As it is mentioned in the paper by \textcite{GortlerHealyThurston2010},
the randomized linear algebra computations,
that at the moment are implemented using integer numbers,
would benefit from a modular approach:
one first solves the equations modulo a prime $p$,
then the obtained solution is lifted to one modulo $p^2$, and so on.
This approach may be considered in some future versions of \pyrigi{}.

\minisec{Extension Constructions}

Common means of constructing rigid graphs are given by the coning operation and
$k$-extensions. \textit{Coning} a graph $G=(V,E)$ adds a new vertex $v^* \notin V$
and edges from $v^*$ to each original vertex $v\in V$.
Coning carries over the generic and global $d$-rigidity of a graph to
$\mathbb{R}^{d+1}$ \parencite{ConnellyWhiteley2010}.
The cone graph can be constructed via the method \pythoninline{cone}.

Given a number $k\in \mathbb{N}$ and a vertex $v\notin V$,
a \textit{$d$-dimensional $k$-extension} of $G=(V,E)$ is obtained by removing from $E$ a subset $F \subset E$ of edges with $|F|=k$,
adding the vertex $v$ to $V$ and, subsequently,
adding edges from $v$ to any subset of $V$ with $d+k$ elements that contains all vertices incident to the edges in $F$.
These extensions are important because in many cases they preserve some rigidity properties.
For instance 0- and 1-extensions preserve $d$-rigidity \parencite{TayWhiteley1985}.

\minisec{Rigidity Matroid}

Rigidity theory has always been intimately connected with matroid theory:
many of the concepts that we have been discussing so far in this paper can be defined in terms of matroids.
Matroids are combinatorial objects that generalize both linear independece and cycles in a graph.
For example, one can define a matroid from a $d$-dimensional framework $(G,p)$ by considering
the linear matroid induced by the rows of its rigidity matrix.
By showing that any two generic $d$-dimensional frameworks $(G,p)$ and $(G,q)$ define the same matroid, one defines the \emph{(generic) rigidity $d$-matroid} of the graph $G$.
The bases of the rigidity $d$-matroid of a complete graph are precisely the minimally rigid graphs on the corresponding number of vertices.

\pyrigi{} offers some basic functionality regarding computations in the rigidity matroid.
The methods \pythoninline{is_Rd_dependent}, \pythoninline{is_Rd_independent} and \pythoninline{is_Rd_circuit} allow the user to determine whether a set of edges is dependent, independent or a circuit in the rigidity matroid.
Similarly, \pythoninline{is_Rd_closed} determines whether a set is closed, and \pythoninline{Rd_closure} computes the closure of a given set in the rigidity matroid.

\minisec{NAC-colorings}
A generically rigid graph might still have non-generic realizations that are flexible.
For any connected graph such a flexible realization in $\RR^2$ exists
if and only if there exists a NAC-coloring \parencite{GraseggerLegerskySchicho2019}.
This is a coloring of the edges in two colors, say red and blue, such that both colors occur
and every cycle is either monochromatic or it contains both colors at least twice.

In general, the existence problem for NAC-colorings is NP-complete \parencite{Garamvölgyi2022}, even for graphs with degree at most five \parencite{LastovickaLegersky2024}.
Nevertheless, the method \pythoninline{NAC_colorings} is able to compute all NAC-colorings for graphs that are not too large.
This is done by determining classes of edges that must be monochromatic in every NAC-coloring.
The naive approach is then to color all these classes in all possible ways and to check which yield a NAC coloring.
In \textcite{LastovickaLegersky2024}, this is improved by an incremental approach considering coloring of subgraphs, supported by new heuristics.
The code supporting \textcite{LastovickaLegersky2024},
which is much faster than the one from \textcite{FlexrilogPaper},
has been included in \pyrigi.

A NAC-coloring is guaranteed to exist if a graph has a stable separating set,
i.e., a set of vertices that have no edges in between and whose removal disconnects the graph.
Such a set can be found by a polynomial time algorithm for 2-flexible graphs \parencite{ClinchGaramvölgyiEtAl2024},
which is implemented in \pyrigi\ within the method \pythoninline{stable_separating_set}.

\subsection{Motions}

Beyond visualizing static frameworks, it is possible to visualize deformation paths of frameworks in \pyrigi.
Given a framework $F=(G,p)$, with $G=(V,E)$ and $p:V\rightarrow \mathbb{R}^d$, a \emph{motion} is a continuous map $\alpha:\mathcal{I}\rightarrow(\mathbb{R}^d)^V$ for an interval $\mathcal{I}\subset \mathbb{R}$ with $0\in \mathcal{I}$ such that $\alpha(0)=p$ and such that $(G,p)$ and $(G,\alpha(t))$ are equivalent for every $t\in \mathcal{I}$. A motion is called \emph{trivial}, if $(G,p)$ and $(G,\alpha(t))$ are congruent for all $t\in \mathcal{I}$.

\minisec{Parametrized Motions}

The most straightforward way to define a continuous motion is to actually provide a parametrization for it.
When constructing a \pythoninline{ParametricMotion} object, we need to specify a parametrization in a single parameter and a (potentially unbounded) interval.
For example, almost all frameworks $F$ on the 4-cycle have a 1-dimensional deformation space that contains the realizations of $F$.
In this case, the corresponding system of polynomials is sufficiently simple so that
one can parametrize a curve in this space by symbolically solving the bar-length equations. At the moment, the specific parametrization of the curve has to be be provided by the user.

\minisec{Numerical Approximations}

However, the system of polynomial equations describing a framework's configuration space quickly becomes too complicated to be described symbolically.
If we still want to visualize a framework's deformation path, we need to rely on numerical methods. Using homotopy continuation enables us to numerically track solution paths \parencite{BreidingTimme2018}. Intuitively, paths are approximated by iteratively applying a \emph{predictor-step} that anticipates the next value on the curve (e.g.\ Euler's method) followed by a sequence of \emph{corrector-steps} to refine the previous prediction, usually by applying Newton's method.

Since the framework's deformation space naturally lives in $(\mathbb{R}^d)^V$,
to approximate a curve in it we can apply the Euclidean distance retraction \parencite{HeatonHimmelmann2025}.
This retraction is given by the metric projection to the closest point on the configuration space.
It works as follows: Given a flexible realization $p$ and a nontrivial infinitesimal flex $q$,
we apply the path-tracking method to the closest point using the Lagrange multiplier system for the parameter $p+\varepsilon q$
and the step size $\varepsilon>0$.
This approach comes with theoretical convergence guarantees, as long as we stay sufficiently close to the constraint set.
There are further theoretical considerations about singularities,
vector transport and stressed frameworks that are incorporated in this method,
though their specifics are beyond the scope of this paper.
In the following example, we approximate a continuous motion of the complete bipartite graph $K_{2,4}$ with $348$ steps,
a step size of $0.1$ and starting with the first infinitesimal flex (\pythoninline{chosen_flex=0})
from a basis that is internally computed so that the pair of vertices $\{0,1\}$ is constrained to the direction $(1,0)$:

\begin{python}
from pyrigi.motion import ApproximateMotion
F = frameworkDB.CompleteBipartite(2,4)
motion = ApproximateMotion(
		F, 348, chosen_flex=0, step_size=0.1, fixed_pair=[0,1]
	  )
\end{python}
The number and size of the steps is chosen so that the approximated continuous motion becomes periodic.

\minisec{Animations}

The computed motions can be animated using the internal \pythoninline{animate} procedure in \pyrigi.
By default, this method produces an SVG animation (restricted to frameworks in $\mathbb{R}$ and~$\mathbb{R}^2$),
but it can also produce a \pythoninline{matplotlib.FuncAnimation} by setting the parameter
\pythoninline{animation_format} to \pythoninline{"matplotlib"}.
This method works regardless of whether it is called on an instance of the \pythoninline{ParametricMotion}
or \pythoninline{ApproximateMotion} class.
While for objects of the former class, the \pythoninline{animate} method needs
to first sample realizations from the parametric curve,
objects from the latter class already come with that information.
Most \pythoninline{PlotStyle} parameters that are admissible for framework plotting can be used here to change the style.

\section{Implementation Details}
\label{implementation}

In this section, we briefly discuss how \pyrigi\ is structured and what external tools it uses.
\pyrigi{} uses several other Python packages:
in particular, we take advantage of the graph theory functionality implemented
in \emph{NetworkX}~\parencite{networkx}. For exact and numerical computations with frameworks
we use \emph{SymPy}~\parencite{sympy} and \emph{NumPy}~\parencite{numpy} respectively,
while \emph{Matplotlib}~\parencite{matplotlib} is used for visualizations.

From a user's perspective, \pyrigi{} follows an object-oriented design.
The functionality of the three areas depicted in \Cref{fig:pyrigi-overview}
can be accessed via the methods of the following classes:
\begin{itemize}
	\item \pythoninline{Graph} --- graph-related functionality, inherited from \pythoninline{networkx.Graph},
	\item \pythoninline{Framework} --- framework-related functionality, inherited from \pythoninline{FrameworkBase}
		(which is composed of \pythoninline{pyrigi.Graph} and a dictionary representing a realization),
	\item \pythoninline{ParametricMotion} and \pythoninline{ApproximateMotion} --- functionality related to continuous motions,
		both inherited from \pythoninline{Motion}.
\end{itemize}
Moreover, the class \pythoninline{GraphDrawer} provides functionality for drawing graphs and 2D frameworks,
and \pythoninline{PlotStyle2D}, \pythoninline{PlotStyle3D} and their parent \pythoninline{PlotStyle}
can be used for specifying style of plots.

\begin{figure}[ht]
	\begin{tikzpicture}
		\node (t) at (0,0) {\begin{minipage}{12cm}
		                  	\dirtree{%
										.1 pyrigi/.
										.2   graph/.
										.3     \_constructions/.
										.4       extensions.py \DTcomment{$k$-extensions}.
										.3     \_rigidity/.
										.4       generic.py \DTcomment{functions for generic rigidity}.
										.4       global\_.py \DTcomment{functions for global rigidity}.
										.4       matroidal.py \DTcomment{functions for rigidity matroid}.
										.4       redundant.py \DTcomment{functions for redundant rigidity}.
										.3     graph.py \DTcomment{implementation of class \pythoninline{Graph}}.
									}
		                  \end{minipage}};
		\draw[con] ($(t.south east)+(0,0.4)$) -| ++(0.6,1) |- node[pos=0.1,below,labelsty,rotate=90] {methods wrap functions from here} ($(t.north east)+(0,-1.8)$);
		\foreach \y in {-2.75,-3.23,...,-4.5}
		{
			\draw[con] ($(t.south east)+(0,0.4)$) -| ++(0.6,0.25) |- ($(t.north east)+(0,\y)$);
		}
	\end{tikzpicture}
	\caption{An example how the package is structured into modules according to topics.}
	\label{fig:module_structure}
\end{figure}
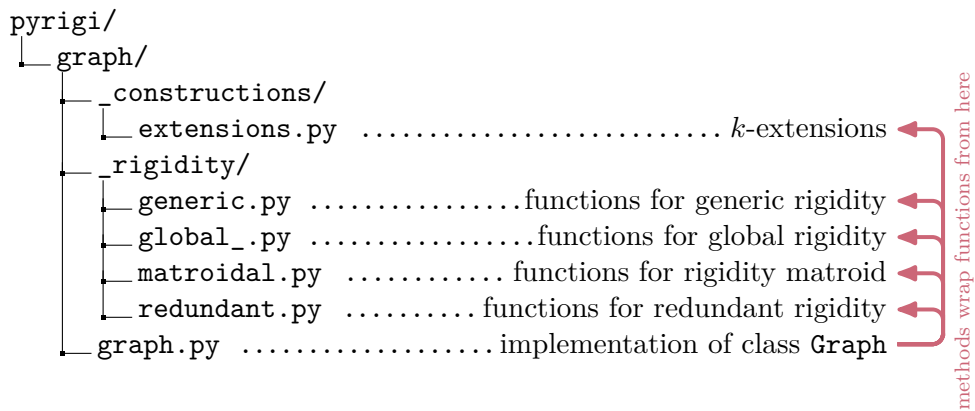

However, in order to have extensible and maintainable code,
most of the graph algorithms are implemented as functions
accepting a \pythoninline{networkx.Graph} instance as the first parameter.
These functions are then wrapped as \pythoninline{pyrigi.Graph} methods.
Hence, they are suggested by autocompletion tools once
a \pythoninline{pyrigi.Graph} instance is available and therefore easy to search for and use.
On the other hand, this approach allows the implementation of functionality in separate non-public modules
according to their topics
which roughly correspond to the structure of \Cref{fig:pyrigi-overview} and \Cref{functionality},
see \Cref{fig:module_structure} for an illustration.
Similarly, the functions operating with frameworks accept \pythoninline{FrameworkBase}
(which implements only the basic methods like adding vertices or edges) and
are then wrapped as methods of \pythoninline{Framework}.

\minisec{Development Tools}

Since \pyrigi{} aims to be a larger longer-term project
that allows incorporating new areas of rigidity theory,
several development tools and processes have been set up.
In this document we just summarize them. More details can be found
in the development guide in the documentation.

Following the current standard in software development, the source code control is done using \emph{Git},
with the remote repository hosted on \emph{GitHub}\footnote{\url{https://github.com/PyRigi/PyRigi}}.
Since libraries are explicitly versioned rather than continuously developed,
the \emph{Gitflow}\footnote{Vincent Driessen, \url{https://nvie.com/posts/a-successful-git-branching-model}} approach is employed, namely, the branch \pythoninline{main}
contains stable code, while the development is done on the branch \pythoninline{dev}
via feature branches.
Whenever appropriate, a new version is released by tagging that
follows a \emph{major.minor.patch} numbering scheme.
Tag creation automatically triggers a GitHub action that deploys the release
to the \emph{Python Package Index}\footnote{\url{https://pypi.org/project/pyrigi}}.
Hence, users can easily install the latest released version of \pyrigi\ via a \pythoninline{pip} command.

Consistent code formatting is achieved by the code formatter \emph{Black}\footnote{\url{https://black.readthedocs.io}},
and compliance with the \emph{PEP8} standard\footnote{\url{https://peps.python.org/pep-0008}} is checked by \emph{Flake8}\footnote{\url{https://flake8.pycqa.org}}.
In order to maintain dependencies, we use \emph{Poetry}\footnote{\url{https://python-poetry.org}},
which also allows easy creation of virtual environments with clean installation of the required dependencies.

We aim to have high test coverage of the code.
For this purpose we distinguish three levels of tests:
the standard ones are checked on every pull request to the development branch (automatically via a GitHub action),
more tests are launched when merging to the branch \pythoninline{main}
in order to create a new stable version (also automatic), and the last group of tests are long ones
that are launched manually to test extensively a specific functionality during the development.
The testing is done using \emph{pytest}~\parencite{pytest}.

We want to create both the math and code documentation as illustrated in \Cref{code_math} simultaneously with a single compilation.
The standard way of documenting Python modules is to use docstrings directly in the source code.
These are compiled by \emph{Sphinx}\footnote{\url{https://sphinx-doc.org}} to an HTML based documentation.
It includes compilation of the mathematical and non-interface parts which are
written in \emph{MyST} format\footnote{\url{https://myst-parser.readthedocs.io/}} with tutorials being generated
from \emph{Jupyter} notebooks via \emph{Myst-NB}\footnote{\url{https://myst-nb.readthedocs.io}}.
The on-line documentation is updated automatically on every merge to the branch \pythoninline{main}.
Notice that function docstrings are placed in corresponding modules,
while the compiled documentation presents the methods wrapping them as this is the interface that we provide to users.
To avoid duplication, the docstring of a method that wraps a function
is copied automatically using a decorator from the wrapped function.
Namely, the actual source code is not as in \Cref{fig:pyrigi-graph,fig:pyrigi-framework},
where the docstrings are for illustration purposes hard copied from the wrapped functions.

\section{Concluding Remarks}
\label{conclusions}

While \pyrigi\ proved to be useful already,
still many important features are not yet present and await implementation.
Among them, we may cite:
\begin{itemize}
 \item gain graphs and support for rigidity respecting symmetry and periodicity;
 \item implementations of different rigidity notions, other than bar-and-joint, such as for instance linearly constrained;
 \item improvements of the efficiency of the currently implemented methods;
 \item implementations of more general criteria for prestress and second-order rigidity.
\end{itemize}
The best way to show interest in the development of \pyrigi{} is to join its Zulip channel\footnote{\url{https://pyrigi.zulipchat.com}} (or contact the maintainers) and propose new features or bug issues there.

Our long-term goal is that \pyrigi{} becomes an environment in which researchers create modules devoted to their specific research needs.
Ideally, these modules would then become part of \pyrigi{} to be easily available to the community.

\addcontentsline{toc}{section}{Acknowledgments}
\section*{Acknowledgments}

We would like to thank Mafalda Dal Cin for her code contribution related to global rigidity,
Hakan Güler for the graph drawing functionality,
Petr Laštovička for NAC-colorings and help with the package structure, and
András Mihálykó for implementing pebble game algorithms.
We would like to acknowledge the contributions to the \pyrigi{} code by
Jose Capco, Angelo El Saliby, Markéta Hošmánková, Daniel Huczala,
Lenka Kušnírová, Jan Rubeš and Johannes Siegele.

We appreciate the comments, tests and contributions to the math documentation by
Sean Dewar, Alison La Porta, Rebecca Monks and Anthony Nixon.

We acknowledge everybody who has contributed to \pyrigi{} in any sense,
including the many discussions we have had,
especially the participants of the workshop Code of Rigidity (March 11--15, 2024),
which was part of the Special Semester on Rigidity and Flexibility at RICAM, Linz, Austria.
We thank RICAM for funding this event.

This research was funded in part by the Austrian Science Fund (FWF) 10.55776/I6233.
For open access purposes, the authors have applied a CC BY public copyright license to any author accepted manuscript version arising from this submission.
M.\,H.\ acknowledges the support by the National Science Foundation under Grant No.\ DMS-1929284 while
the author was in residence at the Institute for Computational and Experimental Research in Mathematics
in Providence, RI, during the Geometry of Materials, Packings and Rigid Frameworks semester program.
J.\,L.\ was supported by the Czech Science Foundation (GAČR), project No.\ 22-04381L.
We would like to thank Zulip\footnote{\url{https://zulip.com}} for the Zulip Cloud Standard plan provided free of charge.

\printbibliography

\end{document}